\documentclass[12pt,]{amsart}
\usepackage{amsmath}
\usepackage{amscd}
\usepackage{amssymb}
\usepackage{latexsym}

\newtheorem{theorem}{Theorem}[subsection]
\newtheorem{lemma}[theorem]{Lemma}
\newtheorem{corollary}[theorem]{Corollary}
\newtheorem{proposition}[theorem]{Proposition}

\newtheorem{remark}[theorem]{Remark}
\newtheorem{definition}[theorem]{Definition}

\newtheorem{slemma}{Lemma}[section]        
\newtheorem{stheorem}[slemma]{Theorem}        
\newcommand{\bgl}{\begin{equation}}         
\newcommand{\egl}{\end{equation}}
\newcommand{\bgln}{\begin{eqnarray}}        
\newcommand{\egln}{\end{eqnarray}}
\newcommand{\bglnoz}{\begin{eqnarray*}}     
\newcommand{\eglnoz}{\end{eqnarray*}}
\newcommand{\btheo}{\begin{theorem}}
\newcommand{\etheo}{\end{theorem}}
\newcommand{\blemma}{\begin{lemma}}
\newcommand{\elemma}{\end{lemma}}
\newcommand{\bproof}{\begin{proof}}
\newcommand{\eproof}{\end{proof}}
\newcommand{\bbew}{\begin{beweis}}
\newcommand{\ebew}{\end{beweis}}
\newcommand{\bremark}{\begin{remark}\em}
\newcommand{\eremark}{\end{remark}}
\newcommand{\bdefin}{\begin{definition}}
\newcommand{\edefin}{\end{definition}}
\newcommand{\bprop}{\begin{proposition}}
\newcommand{\eprop}{\end{proposition}}
\newcommand{\bcor}{\begin{corollary}}
\newcommand{\ecor}{\end{corollary}}
\newcommand{\mn}{\par\medskip\noindent}

\newcommand{\mc}{\mathcal}

\newcommand{\cB}{\mathcal B}
\newcommand{\cC}{\mathcal C}
\newcommand{\cD}{\mathcal D}

\newcommand{\cH}{\mathcal H}

\newcommand{\cJ}{\mathcal J}
\newcommand{\cK}{\mathcal K}
\newcommand{\cL}{\mathcal L}

\newcommand{\sm}{\bf s}
\newcommand{\lori}{\longrightarrow}

\newcommand{\ve}{\varepsilon}

\def\SEMI{\mbox{$\times\kern-2pt\vrule height5pt width.6pt \kern3pt $}}

\newcommand{\halb}{\frac{1}{2}}

\newcommand{\Tr}{{\rm Tr\,}}

\newcommand{\id}{{\rm id}}
\newcommand{\alg}{{\rm alg}}

\newcommand{\Ad}{{\rm Ad\,}}
\newcommand{\lca}{{\rm LCA}}
\newcommand{\ab}{{\rm Ab}}

\renewcommand{\ker}{{\rm ker}\,}
\newcommand{\coker}{{\rm coker}\,}
\def\Cz{\mathbb{C}}

\def\Nz{\mathbb{N}}

\def\Rz{\mathbb{R}}
\def\Zz{\mathbb{Z}}

\newcommand{\ev}{\operatorname{ev}}
\newcommand{\defeq}{\mathrel{:=}}     
\newcommand{\hot}{\mathbin{\hat{\otimes}}}

\title{Algebraic $K$-theory and locally convex algebras}
\author{Joachim Cuntz}
\author{Andreas Thom}
\address{Joachim Cuntz, Mathematisches Institut, Einsteinstr.62, 48149 M\"unster, Germany}
\email{cuntz@math.uni-muenster.de}
\urladdr{http://www.math.uni-muenster.de/u/cuntz/cuntz}
\address{Andreas Thom, SFB 478, Hittorfstr.27, 48149 M\"unster, Germany}
\email{thoman@math.uni-muenster.de}
\urladdr{http://www.math.uni-muenster.de/u/thoman}

\begin{document}

\maketitle \tableofcontents
\newpage
\section{Introduction}
The category of locally convex algebras is very vast. Besides
natural Fr\'{e}chet-algebras arising from differential
calculus on manifolds it also contains algebras defined
purely algebraically (any complex algebra with a countable
basis is a locally convex algebra in a natural way with the
fine topology). The category is therefore very flexible and
allows constructions mixing analytic and algebraic
structures. It seems to cover nearly all examples of (usually
noncommutative) algebras associated with differential
geometric objects, such as algebras of differential or
pseudodifferential operators, algebras of differential forms,
deformation algebras etc., etc. In connection with index
theory and noncommutative geometry it therefore seems highly
desirable to have at our disposal tools such as a computable
version of topological $K$-theory - preferably in bivariant
form as it had been developed for $C^*$-algebras by Kasparov
and others - and associated invariants, also in connection
with the
various classes of ``infinitesimals'' used by Connes \cite{Con}.\\
In\cite{CuDoc}, the first named author had developed such a
theory for the smaller category consisting of all locally
convex algebras that can be topologized by a family of
submultiplicative seminorms (i.e. the category of projective
limits of Banach algebras). This theory has all the desirable
properties and covers already many of the important examples.
However, there are still very natural examples of locally
convex algebras that do not fit into this category.
Particularly relevant examples are given by algebras of
differential operators and in particular, in the simplest
case, by the so called Weyl algebra.\\ In \cite{CuK} a
bivariant $K$-theory $kk^{\alg}$ on the category of all
locally convex algebras was then indeed developed. This
bivariant homology theory again has very good formal
properties, in particular the usual properties of
(differentiable) homotopy invariance, long exact sequences
associated with extensions and stability under tensoring by
the algebra $\cK$ of rapidly decreasing matrices. Therefore
the $kk^{\alg}$-invariants can be computed for many
interesting examples of locally convex algebras as modules
over the graded coefficient ring $R_*=kk_*^{\alg}(\Cz,\Cz)$.
However, the determination of $R_*$ and more generally, a
$K$-theoretic description of $kk_0^{\alg}(\Cz,A)$ remained
open in \cite{CuK} (it was shown though that $R_0$
admits a non-trivial unital homomorphism into $\Cz$).\\
In the present paper we show that the problem of determining
the coefficient ring can be overcome simply by stabilizing by
the Schatten ideal $\cL_p$ rather than by $\cK$. In fact for
the $\cL_p$-stabilized theory $kk^{\cL_p}$ we show that
$kk_0^{\cL_p}(\Cz,A)=K_0(A \hot \cL_p)$ for $p>1$ and thus, in
particular, that $kk_*^{\cL_p}(\Cz,\Cz)=\Zz[u,u^{-1}]$, with a
generator $u$ in degree $2$ (see Corollary \ref{coeff}).
Stabilization by $\cL_p$ is a very natural operation in
Connes' setting for noncommutative geometry. Similar results
hold for stabilization with more general symmetrically normed
multiplicative Banach ideals (see definition
\ref{operatorideal}) including the dual $\cL_{1^+}$ of the
Macaev ideal.
\\ The key to our result is a locally convex version
of an old homotopy invariance theorem proved by Higson in
\cite{Hig} on the basis of previous work by Kasparov and one
of us. One ingredient in our proof of this result is the
development of a setup in which the technology of abstract
Kasparov modules can be applied algebraically
(and thus in particular to locally convex algebras).
\\
Technically, we prove smooth homotopy-invariance of
$M_2$-stable and split-exact functors (see definition
\ref{m2stable}), i.e. in particular of negative algebraic
$K$-theory, for ``weakly $\cJ$-stable'' locally convex
algebras, $\cJ$ being a harmonic Banach ideal (see definition
\ref{har}). All algebras of the form $B \hot \cJ$, where $B$
is an arbitrary locally convex algebra, are weakly
$\cJ$-stable. The result implies that, for such an algebra
$A$, the canonical evaluation maps from the algebra $A[0,1]$
of smooth $A$-valued functions on $[0,1]$ induce the same map
$K_{i}(A[0,1]) \stackrel{ev_0,ev_1}{\to} K_{i}(A)$ for each
non-positive $i \in \Zz$. In particular, since $A[t] \subset
A[0,1]$, we conclude that $A$ is $K_0$-regular. We mention
that the original version of Higson's theorem has been used
in the category of $C^*$-algebras in a similar spirit by
Suslin-Wodzicki
in their proof of the Karoubi conjecture \cite{SW}, see also \cite{Ros}
for a nice exposition of this result and of related ideas.\\
Our conclusions are not obvious, even in the case $i=0$, since
we do not make any assumption concerning openness of the
group of invertibles, stability under functional calculus
etc.. Indeed, as it turns out, it is wrong for $i=1$,
contradicting a claim in \cite{Kar} concerning $K_1(\cL_p)$.
We include a computation of $K_1(\cJ)$ for any harmonic Banach ideal.\\
For higher algebraic $K$-theory, smooth homotopy invariance
fails in general for weakly $\cJ$-stable algebras. This fact
is, from our point of view, due to the lack of excision.
Corti\~{n}as (see \cite{Cor}) has given a long exact
sequence relating the algebraic and topological $K$-theory of such algebras.\\
Our results imply the existence of a homomorphism of graded
rings $kk_*^{alg}(\Cz,\Cz)=R_* \to \Zz[u,u^{-1}]$. It seems
that a proof of the assertion that this homomorphism is an
isomorphism could be derived from the results claimed in
\cite{Tap}. However,
we were not able to verify the arguments in \cite{Tap}.\\
There is independent partially published work by Wodzicki.
There are hints to a deep relationship between the results
announced in \cite{WECM}, see also \cite{Wetal} pp.3, and
implications of our results. It should be interesting to
pursue this further. \mn The second author wants to thank G.
Corti\~{n}as for fruitful discussions about the last section.
We thank Chr. Valqui and G. Corti\~{n}as for helpful comments.
\section{$\cJ$-stable bivariant $K$-theory}
\subsection{Definitions}
By a locally convex algebra we mean an algebra over $\Cz$
equipped with a complete locally convex topology such
that the multiplication $A\times A\to A$ is (jointly) continuous.\\
This means that, for every continuous seminorm $\alpha$ on
$A$, there is another continuous seminorm $\alpha'$ such that
$$\alpha(xy) \leq \alpha'(x)\alpha'(y)$$
for all $x,y \in A$. Equivalently, the multiplication map
induces a continuous linear map $A \hot A \to A$ from the
completed projective tensor product $A\hot A$. All
homomorphisms
between locally convex algebras will be assumed to be continuous.\\
Every Banach algebra or projective limit of Banach algebras
obviously is a locally convex algebra. But so is every
algebra over $\Cz$ with a countable basis if we equip it with
the ``fine'' locally convex topology, see e.g. \cite{CuWeyl}.
The fine topology on a complex vector space $V$ is given by
the family of {\it all} seminorms on $V$.\mn Let $[a,b]$ be an
interval in $\Rz$. We denote by $\Cz[a,b]$ the algebra of
complex-valued $\cC^\infty$-functions $f$ on $[a,b]$, all of
whose derivatives vanish in $a$ and in $b$ (while $f$ itself
may take arbitrary values in $a$ and $b$). Also the
subalgebras $\Cz(a,b], \Cz[a,b)$ and $\Cz(a,b)$ of
$\Cz[a,b]$, which, by definition consist of functions $f$,
that vanish in $a$, in $b$, or in $a$ and $b$, respectively,
will play an important role. The topology on these algebras
is the usual Fr\'echet topology. \mn Given two complete
locally convex spaces $V$ and $W$, we denote by $V\hot W$
their completed projective tensor product (see \cite{T},
\cite{CuWeyl}). We note that $\Cz[a,b]$ is nuclear in the
sense of Grothendieck \cite{T} and that, for any complete
locally convex space $V$, the space $\Cz[a,b]\hat{\otimes}V$
is isomorphic to the space of $\cC^\infty$-functions on
$[a,b]$ with values in $V$, whose derivatives vanish in both
endpoints, \cite{T}, \S\, 51. \mn Given a locally convex
algebra $A $, we write $A[a,b]$, $A [a,b)$ and $ A (a,b)$ for
the locally convex algebras $A \hat{\otimes} \Cz[a,b]$, $A
\hat{\otimes}\Cz[a,b)$ and $A\hat{\otimes}\Cz(a,b)$ (their
elements are $A$ - valued $\mathcal C ^\infty$-functions
whose derivatives vanish at the endpoints). The algebra
$A(0,1)$ is called the suspension of $A$ and denoted by
$\Sigma A$.\mn We denote by $M_n$ the algebra of $n\times
n$-matrices over $\Cz$ and abbreviate, as usual, $A \otimes
M_n$ by $M_n A$. Given continuous homomomorphisms $\phi_1,
\dots, \phi_n:A \to B$, we denote by $\phi_1 \oplus \dots
\oplus \phi_n: A \to M_n(B)$ the diagonal sum of the
$\phi_i$. \mn We also consider the algebra $$M_\infty =
\lim_{\mathop{\lori}\limits_{k}} M_k(\Cz).$$ It is a locally
convex algebra with the fine topology (which is also the
inductive limit topology in the representation as inductive
limit of the algebras $M_k(\Cz)$). Note that, if $A$ is a
finitely generated algebra, every homomorphism from $A$ to
$M_{\infty}$ factors through $M_k$ for some $k\in \Nz$.
Finally, given a locally convex algebra $A$, we denote by
$A^+$ its unitization (as a locally convex space $A^+$ is the
direct sum $A\oplus \Cz$).

\subsection{Banach Ideals}
\label{opideals} Let $H$ be an infinite dimensional separable
Hilbert space. Denote by $\cB(H)$ the algebra of bounded
operators on $H$. Given $p\in [1,\infty)$, the Schatten ideal
$\cL_p \subset \cB(H)$ is defined as
$$ \cL_p = \{ x \in \cB (H) \,|\; \Tr |x|^p < \infty \} $$
Equivalently, a compact operator $x$ is in $\cL_p$ if the
sequence $(\mu_n)$ of its singular values is in
$\ell^p(\Nz)$. $\cL_p$ is a Banach algebra with the norm
$$ \| x\|_p = (\Tr |x|^p)^{1/p}.$$ We denote by $\cL_{\infty}$
the algebra of all compact operators on $H$.

\bdefin\label{operatorideal} A symmetrically normed,
multiplicative Banach ideal $\cJ$ assigns to each
infinite-dimensional separable Hilbert space $H$ a normed
ideal $\cJ(H)$ in the algebra $\cB(H)$ of bounded operators
such that
\begin{itemize}
\item
the norm $\|.\|_\cJ$ on $\cJ (H)$ is complete,
\item
$\|ABC\|_\cJ \leq \|A\| \|B\|_\cJ \|C\|$ for all $A,C \in
\cB(H)$ and $B\in \cJ(H)$,
\item
there is a natural continuous homomorphism $\cJ(H)\hot \cJ(H)
\to \cJ(H\otimes H)$, compatible with the homomorphism
$\cB(H) \hot \cB(H) \to \cB(H \otimes H)$, and
\item
the assignment $H \mapsto \cJ(H)$ is functorial under unitary
transformations.
\end{itemize}
By abuse of notation, in the following we also denote by $\cJ$
the algebra $\cJ(H)$ (which depends on $H$ only up to unitary
isomorphism). \edefin We will abbreviate the term
symmetrically normed, multiplicative Banach ideal to Banach
ideal, since we are considering only such. The main examples
we have in mind are of course the Schatten ideals $\cL_p$ and
the ideal $\cL_{\infty}$, but also the Macaev ideal and its
dual, which we denote by $\cL_{1^+}$. The Banach ideal
$\cL_{1^+}$ is of major importance in Connes' non-commutative
geometry. It consists of all compact operators with singular
values $(\mu_n)$, satisfying $\sum_{i=1}^N \mu_n = O(log(N))$
(see \cite{Con}, pp.439). Although our main result does not
apply to $\cL_1$ it still applies to $\cL_{1^+}$. \mn The
only non-trivial Banach ideal which is closed in operator
norm in $\cB (H)$, is the ideal $\cL_{\infty}$. For our
purposes, in particular for the proof of theorem
\ref{maintheo}, we need to consider a Banach ideal as a
closed ideal in some algebra related to $\cB(H)$. This can be
achieved using the following construction.\\
Given a Banach ideal $\cJ$, we construct the algebra $\cB(H)
\ltimes \cJ(H)$. Topologically, this algebra is defined to be
$\cB(H) \oplus \cJ(H)$ with $l_1$-norm and the multiplication
is such that as many summands as possible are put in
$\cJ(H)$. Since $\cJ(H)$ is symmetrically normed, $\cB(H)
\ltimes \cJ(H)$ is clearly a Banach algebra. Note that
$\cJ(H)$ is a closed ideal in $\cB(H) \ltimes \cJ(H)$ and
that there is a split extension (see definition
\ref{m2stable})
$$ 0 \to \cJ(H) \to \cB(H)\ltimes \cJ(H) \to \cB(H) \to 0.$$
\bdefin \label{har} A Banach ideal $\cJ$ is called harmonic
if $\cJ(H)$ contains a compact operator with singular values
given by the harmonic series for one (and hence each) Hilbert
space $H$. If the Hilbert $H$ space has a standard basis
$\{e_i,i\in \Nz\}$, we denote by $\omega$ the operator which
multiplies the $n$-th basis vector $e_n$ by $1/n$. $\cJ$ is
harmonic if and only if $\omega \in \cJ(l^2(\Nz))$. \edefin
Note that $\cL_p$ for $p>1$ and $\cL_{1^+}$ are harmonic,
whereas $\cL_1$ is clearly not. \mn Let $\cJ$ be a harmonic
Banach ideal. We have the following chain of inclusions of
sub-algebras of $\cB(l^2(\Nz))$
$$M_{\infty} \subset \cK \subset \cL_{1} \subset \cJ \subset \cL_{\infty}$$
in which only the three right-most algebras are Banach
ideals. Clearly, every smooth compact operator is trace
class. The fact, that $\cL_{1}$ is contained in any harmonic
Banach ideal is due to the following chain of inequalities
$$\|A\|_{\cL_1} \geq \sum_{i=1}^k \mu_i \geq k \mu_k.$$ This
implies, given any $A \in \cL_1$, that there is a bounded
operator $B\in \cB(l^2(\Nz))$ such that $U B \omega U^*=
|A|$, for some partial isometry $U$, hence $A \in \cJ$.

\subsection{Bivariant $kk$-theory}

The bivariant homology theory $kk^{\alg}$ associates with any
pair $A,B$ of locally convex algebras abelian groups
$kk^{\alg}_n(A,B)$, $n\in\Zz$, see section \ref{explicitkk}
for a definition and \cite{CuK} for complete proofs. We list
some important properties of $kk^{\alg}$.
\begin{theorem}\label{sixterm}  (a) Every continuous homomorphism
$\alpha : A\to B$ determines an element $kk(\alpha)$ in
$kk_0^{\alg}(A,B)$. Given two homomorphisms $\alpha$ and
$\beta$, we have $kk(\alpha\circ
\beta)=kk(\beta)kk(\alpha)$.\\
(b) Every extension (see section \ref{splexact} for a
definition)
$$
E:\, 0\to  I\stackrel{i}{\lori}  A\, \mathop{\lori}\limits^q
 B\,\to 0
$$ determines canonically an element $kk(E)$ in
$kk_{-1}^{\alg}(B,I)$. The class of the cone extension
$$0\to \Sigma A\to A(0,1]\to A\to 0$$ is the identity element in
$kk_0^{\alg}(A,A)=kk_{-1}^{\alg}(A,\Sigma A)$.\\
If $$
\begin{array}{cccccccccc}
(E):&0 & \rightarrow & A_{1} & \rightarrow & A_{2} &
\rightarrow
& A_{3} & \rightarrow & 0\\
 && & \;\downarrow {\scriptstyle\alpha } & & \downarrow & &
 \downarrow{\scriptstyle\beta } & &\\
(E'):&0 & \rightarrow & B_{1} & \rightarrow & B_{2} &
\rightarrow & B_3 & \rightarrow & 0
\end{array}
$$ is a morphism of extensions (a commutative diagram where the
rows are extensions), then $kk(E)kk(\alpha)=kk(\beta)kk(E')$.\\
(c) $kk^{\alg}$ satisfies Bott periodicity: $kk^{alg}_{n-2}(A,B)\cong
kk^{alg}_n(A,B)$.\\
(d) Let $ D$ be any locally convex algebra. Every extension
$$
E:\, 0\to  I\stackrel{i}{\lori}  A\, \mathop{\lori}\limits^q
 B\,\to 0
$$
induces exact sequences in $kk^{\rm alg}( D,\cdot\,)$ and
$kk^{\rm alg}(\,\cdot\,,  D)$ of the following form: \bgl
\label{1exact}
\begin{array}{ccccc}
kk^{\rm alg}_0( D,  I) & \stackrel{\,\cdot kk(i)}{\lori} &
kk^{\rm alg}_0( D,
 A) &\stackrel{\,\cdot kk(q)}{\lori} & kk^{\rm alg}_0(
D, B)
\\[3pt]
\uparrow  & & & & \downarrow
\\[2pt]
kk^{\rm alg}_1( D,  B) & \stackrel{\,\cdot
kk(q)}{\longleftarrow} & kk^{\rm alg}_1( D,  A) &
\stackrel{\,\cdot kk(i)}{\longleftarrow} & kk^{\rm alg}_1( D,
I)
\end{array}
\egl and \bgl \label{2exact}
\begin{array}{ccccc}
kk^{\rm alg}_0( I, D) & \stackrel{kk(i)\cdot
\,}{\longleftarrow} & kk^{\rm alg}_0(  A,  D) &
\stackrel{kk(q)\cdot \,}{\longleftarrow} & kk^{\rm alg}_0( B,
D)
\\[3pt]
\downarrow  & & & & \uparrow
\\[2pt]
kk^{\rm alg}_1( B,  D) & \stackrel{kk(q)\cdot\,}{\lori} &
kk^{\rm alg}_1( A,
 D) &\stackrel{kk(i)\cdot\,}{\lori} & kk^{\rm alg}_1(
I,  D)
\end{array}
\egl The vertical arrows in (\ref{1exact}) and (\ref{2exact})
are (up to a sign) given by right and left multiplication,
respectively, by the class $kk(E)$.\\
(e) For each locally convex algebra $D$, there is a
multiplicative transformation $\tau_D : kk^{\alg}_*(A,B) \to
kk^{\alg}_* (A\hot D,B\hot D)$ such that
$\tau_D(kk(\alpha))=kk(\alpha\otimes \id_D)$, for any
homomorphism $\alpha:A\to B$.
\end{theorem}
For the proof see \cite{CuK}.
\begin{definition} Let $ A$ and $ B$ be locally convex
algebras and $\cJ$ be a Banach ideal. We define $$
kk^{\cJ}_n( A ,\,  B\,) = kk^{\alg}_n(A,B\hot \cJ)$$
\end{definition}
A priori, the definition depends on an additional variable,
the Hilbert space. To see the independence, we only need to
know that the action of the unitary group of the Hilbert
space $H$ on $\cJ(H)$ is trivial after applying $H \mapsto
kk^{\alg}(A,B \otimes \cJ(H))$. This is easy to see and a
particular case of Lemma \ref{inner}. \\ Having the
independence of the Hilbert space, we can use the natural map
$\cJ(H)\hot \cJ(H)\to \cJ(H\otimes H)$ and its obvious
associativity property in combination with property
\ref{sixterm} (e), to define an associative product
$kk^{\cJ}_n(A,B)\times kk^{\cJ}_m(B,C)\to
kk^{\cJ}_{n+m}(A,C)$. \\ The properties listed in
\ref{sixterm} extend to the bivariant homology theory
$kk^{\cJ}$, since the projective tensor product is exact on
extensions with continuous split. \mn The following theorem,
which is of independent interest, will show that for the
ideals $\cL_p$ ($p\in [1,\infty)$ or $p=1^+$), the theories
$kk^{\cL_p}$ are all naturally isomorphic.
\begin{theorem}\label{square}
Let $A$ and $B$ be locally convex algebras. We assume that
there is a continuous homomorphism $\alpha: A \to B$ and a
continuous map $\beta: B \hot B \to A$, such that
\begin{itemize}
\item
$\beta \circ (\alpha \hot \alpha): A \hot A \to A$ is the
multiplication on $A$ and
\item
$\alpha \circ \beta: B \hot B \to B$ is the multiplication on
$B$.
\end{itemize}
Under these conditions, $[\alpha] \in kk^{\alg}(A,B)$ is
invertible.
\end{theorem}
\bproof Consider the locally convex vector space $\Sigma A
\oplus B$. Let $t\mapsto \phi_t$ be a smooth homeomorphism
from $[0,1]$ to itself with vanishing derivatives at the
endpoints and satisfying $\phi_0=0$. Endowed with the
multiplication
$$(a_t,b) \cdot (a_t',b') = (a_ta_t' + \beta(\alpha(a_t) \otimes \phi_t b'+
\phi_t b \otimes \alpha(a_t')
+\phi_t b \otimes \phi_t b'),bb')$$ it is a locally convex
algebra which we denote by $L(\alpha,\beta)$. There is a
natural homomorphism $A(0,1] \to L(\alpha,\beta)$ given by
$$A(0,1] \ni a_t \mapsto (a_t-\phi_t a_1,\alpha(a_1)).$$ Similarly, there
is a natural map $L(\alpha,\beta) \to B(0,1]$ given by the
assignment
$$L(\alpha,\beta) \ni (a_t,b) \mapsto \alpha(a_t) +\phi_t b \in B(0,1].$$

Consider the following diagram of extensions.
\[ \begin{array}{ccccccccc}
0 & \to & \Sigma A   & \to & A(0,1]          & \to & A          & \to & 0 \\
  &     & \downarrow &     & \downarrow      &     & \downarrow &     &   \\
0 & \to & \Sigma A   & \to & L(\alpha,\beta) & \to & B          & \to & 0 \\
  &     & \downarrow &     & \downarrow      &     & \downarrow &     &   \\ 
0 & \to & \Sigma B   & \to & B(0,1]          & \to & B          & \to & 0
\end{array} \]

By \ref{sixterm}, (b), the commutativity of the diagram
implies that the class in $kk^{\alg}_0(A,B)$ determined by the
extension in the middle is left and right inverse to
$kk(\alpha) \in kk^{alg}(A,B)$. \eproof \bremark The argument
in the proof of \ref{square} shows at the same time that
$\Sigma\alpha$ induces an isomorphism $E(\Sigma A)\cong
E(\Sigma B)$ for any half-exact (see \ref{dexact}) functor
$E$ on the category of locally convex algebras which is
diffotopy invariant.\eremark

\bcor For $1\leq q<\infty$ or $q=1^+$, the natural map
$kk^{\cL_1}_n( A ,\, B\,) \to kk^{\cL_q}_n( A ,\, B\,)$
defines an isomorphism for all $A,B$ and $n$.\ecor

\bproof Note that the multiplication $\cL_{2p} \hot \cL_{2p}
\to \cL_{2p}$ factors continuously through $\cL_p$ by
H\"olders inequality. By \ref{square} the inclusion map
$\alpha : A \hot \cL_p \to B \hot \cL_q$ induces an
isomorphism in $kk^{\alg}$ whenever $p\leq q\leq 2p$ and, by
iteration, whenever $p\leq q$.\eproof

\section{$M_2$-stable and split-exact functors}
\subsection{Split exactness}
\label{splexact} We consider functors from the category of
locally convex algebras with continuous homomorphisms, which
we denote by $\lca$, to the category of abelian groups, which
we denote by $\ab$. There is a natural forgetful functor from
the category of locally convex algebras to the category of
locally convex vector spaces. A sequence $$0 \to A \to B\to C
\to 0$$ of locally convex algebras is called an extension, if
it is a split-extension in the category of locally convex
vector spaces, i.e. the middle term is a topological direct
sum of kernel and co-kernel. An extension is called
split-extension if there is a continuous splitting $C \to B$
which at the same time is a homomorphism.

\bdefin \label{m2stable} Let $E: \lca \to \ab$ be a functor.
\begin{itemize}
\item[(a)]
The functor $E$ is called split-exact, if any split extension
of locally convex algebras is mapped to a split extension of
abelian groups, i.e. if, for every extension $0\to I \to A \to
B\to 0$ of locally convex algebras with a homomorphism
splitting $B\to A$, the induced sequence $0\to E(I)\to
E(A)\to E(B)\to 0$ is exact (and then automatically also
split).
\item[(b)]
The functor $E$ is called $M_2$-stable, if for any locally
convex algebra $A$, the natural inclusion map $j_2: A
\rightarrow A \hot M_2$ which embeds $A$ in the upper left
corner induces an isomorphism under $E$.
\end{itemize}
\edefin

The most important example of a split-exact and $M_2$-stable
functor is of course the algebraic $K$-theory functor $K_0$,
but there are other examples which are of independent
interest. Corti\~{n}as has heavily used the diffotopy
invariance result (Theorem \ref{maintheo}), which we will
derive for all weakly $\cJ$-stable, split-exact and
$M_2$-stable functors, for the functors $HP_*$, $K^{inf}_*$
and $KH_*$ in order to obtain results about the structure of
higher algebraic $K$-theory of weakly $\cJ$-stable algebras,
see \cite{Cor} and the definitions in section
\ref{diffinvtheorem}. \mn Split exact functors satisfy the
following familiar property, the proof of which is trivial.

\bprop \label{add} Let $E$ be a split exact functor and let
$\phi,\psi : A \to B$ be homomorphisms with $\phi(x) \psi(y)
= \psi(x) \phi(y) =0$ for all $x,y \in A$. We have that the
linear map $\phi + \psi:A \to B$ is a homomorphism and that
$E(\phi+\psi)= E(\phi)+E(\psi)$. \eprop

We will show that every split exact and $M_2$-stable functor
admits a pairing with abstract Kasparov modules. Kasparov
modules are a convenient way of encoding the extended
functoriality of $M_2$-stable and split-exact functors. The
underlying construction is the one of quasihomomorphisms,
which were introduced in \cite{CuKK}. Next, we give a brief
introduction to the basics of quasihomomorphisms.

\subsection{Quasi-homomorphisms}
Let $\alpha$ and $\bar{\alpha}$ be two homomorphisms $A\to D$
between locally convex algebras. Assume that $B$ is a closed
subalgebra of $D$ such that $\alpha (x)-\bar{\alpha}(x)\in B$
and $\alpha(x)B\subset B$, $B\alpha(x)\subset B$ for all
$x\in A$. We call such a pair $(\alpha,\bar{\alpha})$ a
quasihomomorphism from $A$ to $B$ relative to $D$ and denote
it by $(\alpha,\bar{\alpha}): A\to B$.

We will show that $(\alpha,\bar{\alpha})$ induces a
homomorphism $E(\alpha,\bar{\alpha}):\, E(A) \to E(B)$ in the
following way. Define $\alpha', \bar{\alpha}':A\to A\oplus D$
by $\alpha'(x)=(x,\alpha(x)),
\,\bar{\alpha}'=(x,\bar{\alpha}(x))$ and denote by $D'$ the
subalgebra of $D\oplus A$ generated by all elements
$\alpha'(x),\ x\in A$ and by $0\oplus B$. We obtain an
extension with two splitting homomorphisms $\alpha'$ and
$\bar{\alpha}'$ :
$$ 0 \to B \to D' \to A \to 0$$
where the map $D' \to A$ by definition maps $(x,\alpha (x))$
to $x$ and $(0,b)$ to $0$. The map $E(\alpha,\bar{\alpha})$
is defined to be $E(\alpha') - E(\bar{\alpha}') : E(A)
\rightarrow E(B)\subset E(D')$ (this uses split-exactness).
Note that $E(\alpha,\bar{\alpha})$ is independent of $D$ in
the sense that we can enlarge $D$ without changing
$E(\alpha,\bar{\alpha})$ as long as $B$ maintains the
properties above.

\bprop\label{qhp} The assignment $(\alpha,\bar{\alpha})\to
E(\alpha,\bar{\alpha})$ has the following properties:
\begin{itemize}
\item[(a)] $E(\bar{\alpha},\alpha) = -E(\alpha,\bar{\alpha})$
\item[(b)] If the linear map $\varphi = \alpha - \bar{\alpha}$
is a homomorphism and satisfies $\varphi(x)\bar{\alpha}(y) =
\bar{\alpha}(x)\varphi(y) = 0$ for all $x,y\in A$, then
$E(\alpha,\bar{\alpha}) = E(\varphi)$.
\end{itemize}
\eprop \bproof (a) This is obvious from the definition. (b)
This follows from proposition \ref{add} and the fact that
$\varphi + \bar{\alpha} = \alpha$. \eproof \bdefin\label{dinn}
We say that a functor $E: \lca \to \ab$ is
\begin{itemize}
\item
invariant under inner automorphisms, if $E(\Ad U)=
\id_{E(A)}$ for any invertible element $U$ in a unital algebra
$A$,
\item
invariant under idealizing automorphisms, if the
following holds:\\
Whenever $A$ is isomorphic to a subalgebra of a unital
algebra $D$ and $U$ is an invertible element in $D$ such that
$AU, UA, AU^{-1}, U^{-1}A \subset A$, then $E(\Ad
U)=1_{E(A)}$ for the automorphism $\Ad U: A\to A$.
\end{itemize}
\edefin \blemma\label{inner}
\begin{itemize}
\item[(a)] Every split-exact functor that is invariant under inner
automorphisms is also invariant under idealizing
automorphisms.
\item[(b)] Every $M_2$-stable functor is invariant
under inner automorphisms.
\item[(c)] Every split-exact and $M_2$-stable functor is invariant under
idealizing automorphisms.
\item[(d)] Assume that $E$ is invariant under inner automorphisms and let
$(\alpha,\bar{\alpha}):A\to B$ be a quasihomomorphism
relative to $D$. If there is an invertible element $U$ in $D$
such that $BU, UB, BU^{-1}, U^{-1}B\,\subset \,B$ and such
that moreover $U\alpha (x) - \alpha (x) \,\in\,B$,
$\bar{\alpha}(x)=U\alpha(x)U^{-1}$ for all $x\in A$, then
$E(\alpha,\bar{\alpha})=0$.
\end{itemize}
\elemma \bproof (a) Let $A$ be a locally convex algebra and
$U \in D$ be an idealizing element of $A \subset D$. Let $B$
be the sub-algebra of $D \oplus \Cz[t,t^{-1}]$ which is
generated by $A \oplus 0$ and $U \oplus t$ (we view
$\Cz[t,t^{-1}]$ as a locally convex algebra with the fine
topology). There is a natural split extension
$$0 \to A \to B \to \Cz[t,t^{-1}] \to 0.$$ The inner automorphism
determined by $U \oplus t$ induces the identity on $E(B)$ and
thus on $E(A)$.

(b) Let $U$ induce an inner automorphism of the algebra $A$.
Consider the inner automorphism of $A \hot M_2$ induced by $U
\oplus 1$. Since, by an easy argument, the inclusions of $A$
into the upper left corner and the inclusion into the lower
right corner of $A \hot M_2$ determine the {\it same}
isomorphism under $E$, this implies that $E$ is invariant
under idealizing automorphisms.

(c) Combine (b) and (a).

(d) Note that $W=(1,U)$ defines an idealizing automorphism of
the algebra $D'$ that appears in the split-extension defining
the map $E(\alpha,\bar{\alpha})$. By (c), we know that
$E(\bar{\alpha}')=E(\Ad W) \circ E(\alpha')=E(\alpha')$.
Hence, we get $E(\alpha,\bar{\alpha})=0$. \eproof

\subsection{Kasparov modules}
\bdefin \label{KM}Let $A,I$ and $D$ be locally convex
algebras. Assume that $D$ is unital and contains $I$ as a
closed ideal. An abstract Kasparov $(A,I)$-module relative to
$D$ is a triple $(\varphi,U,P)$ where
\begin{itemize}\item $\varphi$ is a continuous homomorphism
from $A$ into $D$.\item $U$ is an invertible element and $P$
is an idempotent element in $D$ such that the following
commutators are in $I$ for all $x\in A$:$$[U,\varphi (x)],\;
[P,\varphi (x)],\; [U,P].$$\end{itemize}\edefin

\bremark (a) The condition that $I$ is a closed ideal in $D$
could be weakened to demanding that the inclusion $I \to D$
and the homomorphisms $A \otimes I \to I$ and, $I \otimes A
\to I$ which are given by mapping $A$ to $D$ and multiplying
in $D$ are continuous for the projective topology on the
tensor product.

Furthermore, just as for quasihomomorphisms, the morphisms
that we construct from a Kasparov module will be independent
of the algebra $D$ in the sense that $D$ can be enlarged as
long as the conditions on the commutators and the ideal are
maintained.

(b) In Definition \ref{KM} one could replace the condition
that $[U,P]\in I$ by the condition that $[U,P]\varphi (x)$ and
$\varphi (x)[U,P]$ are in $I$ for all $x\in A$. \eremark

In the sequel, given an idempotent $P$ in a unital algebra,
we will denote by $P^\perp$ the idempotent $1-P$. Given an
abstract Kasparov module $(\varphi,U,P)$, we define invertible
matrices $W_P$ and $U_P$ in $M_2(D)$ as follows
$$
 W_P=\left(\begin{array}{cc}
 P & P^\perp\\
 -P^\perp & P \end{array} \right),
$$ and
  $$U_P= W_P\left(\begin{array}{cc}
 U & 0\\  0 & 1 \end{array} \right) W_P^{-1} = \left(\begin{array}{cc}
 PUP+P^\perp & -PUP^\perp\\  -P^\perp UP & P+P^\perp UP^\perp \end{array} \right). $$\mn

We obtain in this way a quasihomomorphism $((\varphi\oplus
 0),{\rm Ad}\, U_P\circ (\varphi\oplus
 0)):\,A\to M_2(I)$.
 \bdefin Let $E$ be a functor that is split exact and
 $M_2$-stable and let $(\varphi,U,P)$ be an abstract
 $(A,I)$-Kasparov module. We denote by $E(\varphi,U,P)$ the
 map $E(A)\to E(M_2 I)\cong E(I)$ associated with the
 quasihomomorphism $((\varphi\oplus 0),{\rm Ad}\, U_P\circ (\varphi\oplus
 0)):\,A\to M_2 I$ and the natural inclusion $j_2:I \to M_2 I$.\edefin
 \blemma\label{modprop} Let $E$ be split exact and $M_2$-stable. We consider
 abstract $(A,I)$-Kasparov modules.
 \begin{itemize}
 \item[(a)] $E(\varphi,UU',P)=E(\varphi,U,P) + E(\varphi,U',P)$ whenever
 all terms are defined.
 \item[(b)] If $PP'=P'P=0$, then $E(\varphi, U,P+P')=
 E(\varphi, U,P)+E(\varphi, U,P')$.
 \item[(c)] If $P$ commutes with all $\varphi (x)$,
 $\varphi'(x)$ and if $P\varphi (x)= P\varphi'(x)$ for all $x$,
 then $E(\varphi,U,P)=E(\varphi',U,P)$. \item[(d)]
 If $\varphi(x)U-\varphi(x)\,\in I$ for all $x\in A$, then $E(\varphi,U,P)=0$.
 \item[(e)] If $PU - P'U' \in I$, then $E(\varphi,U,P)=E(\varphi,U',P')$.
 \end{itemize}

 \elemma \bproof
 (a) We have $(UU')_P=U_PU'_P$. This implies
\[ \begin{array}{rcl}
E(\varphi,UU',P)&=& E(\varphi \oplus 0, Ad((UU')_P) \circ
(\varphi \oplus 0) )
\\ &=&   E(\varphi \oplus 0, Ad(U_P) \circ (\varphi \oplus 0))
\\ & & + \,E(Ad(U_P)) \circ (\varphi \oplus 0),Ad(U_PU'_P) \circ (\varphi \oplus 0))
\\ &=& E(\varphi,U,P) + E(Ad(U_P)) \circ E(\varphi,U',P)
\\ &=& E(\varphi,U,P) + E(\varphi,U',P)
\end{array} \]
The last equality follows, since $Ad(U_P)$ is an idealizing
automorphism of $M_2(I)$ and we have shown the invariance
under such in Lemma \ref{inner}(c).

(b) One easily checks that
 $1-(U_{P+P'})^{-1}U_PU_{P'}\in M_2(I)$. The assertion then
 follows from (a) combined with (d).

(c) A computation using the definition of a Kasparov module
and $M_2$-stability shows that $E(P^\perp\varphi,
U,P)=E(P^{\perp}\varphi \oplus 0, P^\perp\varphi \oplus 0)=0$.

This observation implies that  $$E(\varphi, U,P)=E(P\varphi,
U,P) = E(P\varphi',U,P)=E(\varphi', U,P).$$

(d) The element $U_P$ satisfies the hypotheses of lemma
 \ref{inner}(d) with respect to $\varphi\oplus 0$.

 (e) We
 have the following chain of congruences mod the ideal $I$:
 $$P\,\equiv\,U^{-1}U'P'\,=\,(U^{-1}U'P')P'\,
 \equiv\,PP'\,\equiv\,P(PUU'^{-1})
 \,=\,PUU'^{-1}\,\equiv\,P'$$
 ($\equiv$ denotes congruence mod $I$).\\
 The explicit form of the matrix $U_P$ is as follows:
 $$U_P \;= \left(\begin{array}{cc}
  PUP+P^\perp &- PUP^\perp\\
 -P^\perp UP & P+P^\perp UP^\perp \end{array} \right)$$
 and similarly for $U'_{P'}$. Therefore the matrices
$${U'}^{-1}_{P'} U_P (1 \oplus 0)$$ and
$$(1 \oplus 0) U_P^{-1} U'_{P'}$$ are in $M_2(I^+)$.
These two matrices produce an algebraic equivalence in
$M_2(I^+)$ between the idempotents $1 \oplus 0$ and
${U'}_P^{-1}U_P(1 \oplus 0)U_P^{-1}U'_{P'}$. Using the
standard technique (e.g. \cite{Bla} prop. 4.3.1) there exists
an invertible element $V$ in $M_4(I^+)$ such that
$$V ({U'}_{P'}^{-1} U_P \oplus 1 \oplus 1) (1 \oplus 0 \oplus 0
\oplus 0)  = (1 \oplus 0 \oplus 0 \oplus 0)$$ and
$$(1 \oplus 0 \oplus 0 \oplus 0) (U'_{P'} U_P^{-1} \oplus 1
\oplus 1)V^{-1}  = (1 \oplus 0 \oplus 0 \oplus 0).$$
Moreover, $V$ can be chosen such that $V \equiv 1 \oplus
1\oplus 1\oplus 1$. This last observation, together with
$M_2$-stability, implies the third equality in the following
computation.

\[ \begin{array}{rcl}
0 &=& E(\varphi \oplus 0 \oplus 0 \oplus 0,\varphi \oplus 0 \oplus 0 \oplus 0)\\
  &=& E(\varphi \oplus 0 \oplus 0 \oplus 0, Ad(V) \circ Ad({U'}^{-1}_{P'} U_P \oplus 1 \oplus 1) \circ (\varphi \oplus 0 \oplus 0 \oplus 0))\\
  &=& E(\varphi \oplus 0, Ad({U'}^{-1}_{P'} U_P) \circ (\varphi \oplus 0))\\
  &=& E(\varphi \oplus 0, Ad({U'}^{-1}_{P'})(\varphi \oplus 0)) + E(Ad({U'}^{-1}_{P'})) \circ E(\varphi \oplus 0,Ad(U_P) \circ (\varphi \oplus 0 )) \\
  &=& E(\varphi,{U'}^{-1},P') + E(\varphi,U,P)
\end{array} \]
The last equality is due to the observation that
$Ad({U'}^{-1}_{P'})$ is an idealizing automorphism of
$M_2(I)$. This proves the assertion
$$E(\varphi,U',P')=E(\varphi,U,P),$$ using (a) and (d).
\eproof

\subsection{Diffotopy invariance via Kasparov
modules}\label{difinv}

In this section we construct explicit Kasparov modules which
induce the evaluation maps, stabilized by the smooth compact
operators. They are used in the proof of the diffotopy
invariance theorem (see Theorem \ref{maintheo}).\mn

We let $\cC^\infty (S^1)$ act in the usual way on
$\ell^2(\Zz)$ and denote by $P$ the Hardy projection from
$\ell^2(\Zz)$
onto the subspace $\ell^2(\Nz)$.\\
We denote by $D$ the subalgebra of $\cB(\ell^2(\Zz))$
generated algebraically by $\cC^\infty (S^1)$ together with
$P$ and all smooth compact operators. There is a natural
isomorphism of $\cC^\infty (S^1)$ with the Schwartz space of
rapidly decreasing $\Zz$-sequences $\sm\subset \ell^2(\Zz)$
defined by the action on the basis vector $e_0 \in
\ell^2(\Zz)$. A vector is called smooth if it
belongs to the image of this embedding.\\
We list some elementary properties concerning the action of
the smooth algebras that we consider on $\ell^2(\Zz)$.
\begin{lemma}
\label{smooth-lemma}
\begin{itemize}
\item[(a)]There is an isomorphism of topological vector spaces
$\sm \hot \sm \to \cK$ extending the map $(\lambda_i)_{i\in
\Zz} \otimes (\mu_j)_{j\in\Zz} \mapsto (\lambda_i\mu_j)_{i,j
\in \Zz}$. In particular, any finite-rank operator built out
of smooth vectors is a smooth compact operator.
\item[(b)] The action of $\cK\cong\sm\hot\sm$ on the Hilbert space
$\ell^2(\Zz)$ is induced by $a\otimes b(\xi)=a\langle
b|\xi\rangle$ where $\langle \cdot|\cdot\rangle$ denotes the
canonical bilinear form on $\ell^2(\Zz)$. In particular the
image $k(\ell^2(\Zz))$ of any $k\in\cK$ is contained in
$\sm$. Any eigenvector for an eigenvalue $\neq 0$ of an
element $k\in\cK$ is smooth. \item[(c)] The algebra $\cK$
(resp. its unitization $\cK^+$) is a subalgebra of the
$C^*$-algebra $\cL_\infty$ of compact operators (resp. of its
unitization $\cL_\infty^+$) which is closed under functional
calculus by holomorphic functions for arbitrary elements (in
particular the spectrum of an element in $\cK$ is the same as
in $\cL_\infty$), and under functional calculus by
$\cC^\infty$-functions for self-adjoint and normal elements.
\item[(d)] The algebra $D$ preserves the space of smooth vectors.
\end{itemize}
\end{lemma}
\bproof The only non-trivial assertion is (c). Let $N$ be the
unbounded operator in $H=\ell^2(\Zz)$ defined by
$N(e_i)=|1+i|e_i$ on the domain $\cD$ consisting of all
(finite) linear combinations of vectors $e_i$ in the standard
orthonormal basis for $H$. It is trivial to check that $\cK$
can be identified with the set of operators $a$ in $\cB (H)$
for which $N^i aN^j$ extends from $\cD$ to a bounded operator
on $H$
for all $i,j\in \Nz$.\\
Define norms $\eta_k$ on $\cK$ by
$\eta_k(a)=\sum_{i+j=k}\|N^i a N^j\|$ where $\|\,\cdot\,\|$
denotes the operator norm. Then, for $a,b\in \cK$ we have:
$$\eta_k(ab)=\sum_{i+j=k}\|N^i ab N^j\|\leq  \sum_{i+j=k}\|N^i a\|\,\|b N^j\|
\leq \sum_{i+j=k}\eta_i(a)\eta_j(b)$$ Thus, the family
$(\eta_j)$ defines a differential seminorm in the sense of
\cite{BC}. The topology on $\cK$ is described by the family
of sub-multiplicative norms $\|\,\cdot\,\|_k$ defined by
$$ \| a \|_k = \sum_{0\leq j\leq k}\eta_j (a) $$
The results in \cite{BC} then immediately imply that $\cK$ is
a subalgebra of the algebra $\cL_\infty$ of compact operators
on $H$ which is closed under functional calculus by
holomorphic functions for arbitrary elements (see \cite{BC},
3.12) and by functional calculus by $\cC^\infty$-functions for
self-adjoint or normal elements (see \cite{BC}, 6.4).\eproof
Linearly, $D$ splits into $\cK \oplus \cC^\infty(S^1) \oplus
\cC^\infty(S^1)$ via the linear isomorphism
$$\cK \oplus \cC^\infty(S^1) \oplus \cC^\infty(S^1) \ni (k,f,g)
\mapsto k + f + gP \in D.$$ We topologize $D$ by norms
$\|\,\cdot\,\|_k \oplus p_k \oplus p_k$, where
$\|\,\cdot\,\|_k$ is of the type described in Lemma
\ref{smooth-lemma}, (c) and $p_k$ is the usual
sub-multiplicative norm on $\cC^\infty(S^1)$ defined by
$$ p_k(f) = \sum_{0\leq j\leq k}\frac{1}{j!}\|f^{(j)}\|$$ An
easy computation shows that these norms are
sub-multiplicative. This shows that $D$ is a locally convex
algebra and contains the algebra $\cK$ of smooth
compact operators as a closed ideal.\\
We identify $\cC^\infty(S^1)$ with the algebra of smooth
 periodic functions on the interval $[0,4]$ and divide this
 interval into four succeeding intervals $I_1,I_2, I_3,I_4$ of length
 $1$. We define a homomorphism $\varphi :\Cz[0,1]
 \to\cC^\infty(S^1)$ in the following way:
 \bglnoz  \varphi (f) \equiv f(0) \quad{\rm on}\, I_1\qquad\quad
 \varphi (f) = f\quad {\rm on} \, I_2\\[2mm]
 \varphi (f) \equiv f(1) \quad{\rm on}\, I_3\qquad\quad
 \varphi (f) = \check{f} \quad{\rm on}\, I_4\eglnoz
 where we identify the interval $[0,1]$ with each of the
 intervals $I_k$ and where we put $\check{f}(t)=f(1-t)$, for
 $t\in [0,1]$.\\
 Let now $u_t$, $t\in [0,1]$ be a smooth family of unitary
 elements in $\cC^\infty (S^1)$ such that each $u_t$ has
 winding number 1 and such that $u_0\equiv 1$ outside a closed
 interval $I_1'$ contained in the interior of $I_1$
 and $u_1\equiv 1$ outside an
 interval $I_3'$ contained in the interior of $I_3$ (there is an obvious
 explicit family with these properties).\\
 Since $[P,\varphi (x)]\in \cK$ and $[P,u_t]\in \cK$ for each
 $x$ and $t$, $(\varphi, u_t, P)$ defines an abstract Kasparov
 module for each $t$.
 \blemma \label{kasparov-evaluation-lemma}
Let $E$ be an $M_2$-stable and split exact functor. Then
 $$E(\varphi,u_0,P)=E(j\circ \ev _0)\qquad E(\varphi,u_1,P)=E(j\circ \ev
 _1)$$ where $\ev_t :\Cz [0,1]\to \Cz$ are the evaluation maps
 and $j:\Cz\to\cK$ is the natural inclusion.
 \elemma\bproof Let $h$ be a function in $\cC^\infty(S^1)$
 such that $h\equiv 1$ on $I_1'$, $h\equiv 0$ outside $I_1$ and
 $0\leq h\leq 1$ and put
 $$ Q\,=\,\left(\begin{array}{cc}
 h & (h(1-h))^\halb\\[2mm]
 (h(1-h))^\halb & 1-h \end{array} \right)$$
 We consider $Q$ as an idempotent in $M_2(D)$.
 It commutes exactly with $\varphi (f)\oplus\ev _0 (f),\,f\in\Cz[0,1]$
 and we have
 $$Q(\varphi(f) \oplus \ev_0(f))= Q(\ev_0(f) \oplus \ev_0(f))$$ and
 $$(u_0 \oplus 1) Q^{\perp} = Q^{\perp}(u_0\oplus 1) = Q^{\perp}.$$
We have, using the last identity and Lemma \ref{modprop} (e),
the following equalities:
$$\begin{array}{rcl}
E(\varphi,u_0,P)\quad & = & E(\varphi \oplus \ev_0, u_0 \oplus 1, P \oplus P)\\
&=& E((\varphi \oplus \ev_0)Q, u_0 \oplus 1, P \oplus P)\\
&=& E((\ev_0 \oplus \ev_0)Q, u_0 \oplus 1, P \oplus P) \\
&=& E(\ev_0 \oplus \ev_0,u_0 \oplus 1, P \oplus P) \\
&=& E(\ev_0,u_0,P)
\end{array}$$
An analogous argument shows that $E(\varphi,u_1,P)=E(\ev_1,u_1,P)$.\\
Now, note that $u_0$ is connected to the standard unitary $z=
\{ t \mapsto e^{2\pi i t}\}$ by a continuous path $w_t$ in
$C^{\infty}(S^1)$ with $w_0=z$ and $w_1=u_0$. Since $P$
commutes with all elements in $C^{\infty}(S^1)$ modulo
compact operators we see that $Pw_tP$ describes a continuous
path of operators, which are Fredholm on $\ell^2(\Nz)$. The
operator $PzP$ is clearly of index $-1$, then so is $v=Pu_0P$.
Applying functional calculus by the holomorphic function,
which is equal to $z^{-1/2}$ in a neighborhood of 1 and equal
to 0 in a neighborhood of 0, to $v^*v$, there is a positive
element $h$ in the algebra generated by $P\cK P$ and $P$,
which differs from $P$ only by an element in $\cK$ such that
$vh$ is a partial isometry still of index $-1$ (in fact with
the same kernel and cokernel as $v$). Choose a partial
isometry $k$ in $\cB(\ell^2(\Nz))$ with support $\ker v$ and
range contained in $\coker v$. Since the range and support of
$k$ consist of smooth vectors, $k$ is in $\cK$. Thus $a=vh+k$
is an isometry of index -1 on $\ell^2(\Nz)\subset \ell^2(\Zz)$
which differs from $v$ only by an element in $\cK$.\\
In the same way one gets a coisometry $b$ of index 1 on
$\ell^2(\Zz) \ominus \ell^2(\Nz)$ and equal to $P^\perp u_0
P^\perp$ modulo $\cK$. Adding a partial isometry of rank 1 in
$\cK$, that connects $\ker b$ and $\coker a$, to $a+b$ one
obtains a unitary $u$ in $\cB(\ell^2(\Zz))$ that differs from
$u_0$ only by an element in $\cK$ such that
$uP=a$. \\
Using Lemma \ref{modprop} (e) and the explicit form of the
quasihomomorphism defined by the abstract Kasparov module
$(\ev_0,u,P)$, we have
$$E(\ev_0,u_0,P)=E(\ev_0,u,P)=E(\ev_0,\Ad (a+P^\perp)\ev_0)$$
However, by the choice of $a$, we have $\ev_0 =\Ad
(a+P^\perp)\ev_0 + j\circ\ev_0$, so that $E(\ev_0,u_0,P)=
E(j_0\circ \ev_0)$ by \ref{qhp} (b). The argument for
$E(\ev_1,u_1,P)$ is exactly symmetric.\eproof

\section{Diffotopy invariance theorem} \label{diffinvtheorem}
\subsection{Half-exact functors and diffotopy invariance}

We continue with some definitions.

\bdefin\label{dexact} Let $E: \lca \to \ab$ be a functor from
the category of locally convex algebras to the category of
abelian groups. Let $\cJ$ be a Banach ideal. We say that
\begin{itemize}
\item $E$ is
diffotopy invariant, if the maps $\ev_t: E(A[0,1])\to E(A)$
induced by the different evaluation maps for $t\in [0,1]$ are
all the same (it is easy to see that this is the case if and
only if the map induced by evaluation at $t=0$ is an
isomorphism).
\item $E$ is half-exact, if, for every extension $0\to I \to
A \to B\to 0$ of locally convex algebras, the induced short
sequence $E(I)\to E(A)\to E(B)$ is exact.

\item $E$ is
$\cJ$-stable, if the map $E(A)\to E(A \hot \cJ)$ induced by
the natural inclusion $j:A\to A \hot \cJ$, defined by any
minimal idempotent in $\cJ$, is an isomorphism for each
locally convex algebra $A$ (note that any Banach ideal
contains minimal idempotents).

\item $E$ is weakly $\cJ$-stable, if there is
a natural map $E(A \hot \cJ) \to E(A)$, such that the
composition $E(A)\to E(A \hot \cJ)\to E(A)$ with the map
induced by the inclusion map $A\to A\hot\cJ$ defined by a
minimal idempotent in $\cJ$ is the identity for each locally
convex algebra $A$.
\end{itemize}\edefin

\bremark \label{heredweakstab} We are using weak
$\cJ$-stability only for Banach ideals $\cJ$, but note that
weak stability is hereditary for certain subalgebras of $\cJ$.
In particular, if $\cJ$ contains all smooth compact operators
(for instance if $\cJ$ is harmonic, see end of section
\ref{opideals}), then any weakly $\cJ$-stable functor is also
weakly $\cK$-stable. \eremark

\bdefin\label{dtop} Two homomorphisms $\alpha, \beta :A\to B$
between locally convex algebras are called diffotopic if
there is a homomorphism $\varphi : A \to B[0,1]$ such that
$$ \alpha = \ev_0 \circ \varphi \qquad \beta
=\ev_1 \circ\varphi.$$ \edefin

\bremark In \cite{CuDoc} and elsewhere it was stated that the
condition defining diffotopy in \ref{dtop} was equivalent to
the fact that there exists a family of homomorphisms
$\varphi_t :A\to B,\quad t\in [0,1]$ such that, for each $x$
in $A$, the function $t\mapsto \varphi_t(x)$ is in $B[0,1]$
and such that $\alpha =\varphi_0$, $\beta =\varphi_1$ . The
equivalence between the two conditions holds only if $A$ and
$B$ are Fr\'{e}chet, i.e. metrizable. A counterexample in the
general case, due to Frerick and Shkarin, was communicated to
us by L.Frerick \cite{FS}. Since only the condition given in
the definition above was used in all the proofs, this does not
affect any of the results in \cite{CuDoc} or \cite{CuWeyl}.
\eremark

It is easy to check that diffotopy is an equivalence
relation. If $\alpha$ and $\beta$ are diffotopic and $E$ is a
diffotopy invariant functor, then obviously
$E(\alpha)=E(\beta)$. \mn The next lemma illustrates the
important consequences of the combination of the properties
of diffotopy invariance and half-exactness. Note that the
algebraic $K$-theory functor $K_0$ is half-exact, in fact, it
satisfies excision (see \cite{Bass}). Using the formally
similar property of split-exactness, we are going to conclude
that the $\cJ$-stabilized version of the algebraic $K$-theory
functor is diffotopy invariant. However, split-exactness does
not imply half-exactness in general, nor does the reverse
implication hold. For some of the notation, the reader is
referred to section \ref{explicitkk}.

\blemma \label{halfexact} Let $E: \lca \to \ab$ be a functor
from the category of locally convex algebras to the category
of abelian groups which is diffotopy invariant and
half-exact. Then
\begin{itemize}
\item[(a)] $E$ has long exact sequences, i.e. for each
extension $0\to I \to A \to B\to 0$ of locally convex
algebras there is a long exact sequence (infinite to the
left) of the form
$$ \cdots \to E(\Sigma A)\to E(\Sigma B)\to E(I)\to E(A)\to E(B)$$
\item[(b)] There is a natural isomorphism $E(JA)\cong E(\Sigma A)$
\end{itemize}\elemma \bproof (a) is well known and follows from a standard
argument using
mapping cones.\\
(b) Apply the long exact sequence from (a) to the extension
$0\to JA \to TA \to A\to 0$ and use the fact that
$E(TA)=E(\Sigma TA)=0$ by diffotopy invariance of $E$ and the
fact that $TA$ is smoothly contractible.\eproof

\subsection{Main theorem}

In this section we state and prove the diffotopy invariance
theorem.

\btheo \label{maintheo} Every functor from the category of
locally convex algebras to the category of abelian groups
which is split exact, $M_2$-stable and weakly $\cJ$-stable
for some harmonic Banach ideal, is diffotopy invariant. \etheo

\bproof Given a weakly $\cJ$-stable, $M_2$-stable and
split-exact functor $E: \lca \to \ab$, these properties are
inherited by the functor $E^A: \lca \to \ab$ which assigns $B
\mapsto E(A \hot B)$, for any given locally convex algebra
$A$. This implies that w.l.o.g. it is enough to show that the
two evaluation maps from $E(\Cz[0,1])$ to $E(\Cz)$ are equal.

Using the result of Lemma \ref{kasparov-evaluation-lemma} and
the injectivity of $E(j): E(\Cz) \to E(\cK)$, which followed
from weak $\cJ$-stability and remark \ref{heredweakstab}, it
suffices to show that we have an equality of
$(\Cz[0,1],\cK)$-Kasparov modules
$$E(\varphi,u_0,P)=E(\varphi,u_1,P).$$ Using the weak
$\cJ$-stability again we conclude that it also suffices to
show that the following $(\Cz[0,1],\cK \hot \cJ))$-Kasparov
modules are equal:

$$E(\varphi \otimes 1,u_0 \oplus 1 \oplus \dots, P \otimes 1)
= E(\varphi \otimes 1,u_1  \oplus 1 \oplus \dots, P  \otimes
1).$$

Note that $u_0 {u_1}^{-1} = e^{ih}$ for some self-adjoint
element $h \in C^{\infty}(S^1)$. Denote by $Z$ the unitary in
$(C^{\infty}(S^1) \hot \cJ)^+\subset (D\hot \cJ)^*$ given by
$$Z=e^{ih} \oplus 1 \oplus 1 \oplus \dots.$$
where $D$ is as in \ref{difinv}.

We have to show that $E(\varphi \otimes 1, Z, P\otimes 1)=0$.
The element $Z$ is the product of elements
$$X_1 = e^{ih} \oplus e^{-\frac{ih}{2}} \oplus
e^{-\frac{ih}{2}} \oplus e^{\frac{ih}{4}} \oplus
e^{\frac{ih}{4}} \oplus e^{\frac{ih}{4}} \oplus
e^{\frac{ih}{4}} \oplus e^{-\frac{ih}{8}} \oplus \dots $$ and
$$X_2 = 1 \oplus e^{\frac{ih}{2}} \oplus e^{\frac{ih}{2}} \oplus e^{-\frac{ih}{4}} \oplus e^{-\frac{ih}{4}}
\oplus e^{-\frac{ih}{4}} \oplus e^{-\frac{ih}{4}} \oplus
e^{\frac{ih}{8}} \oplus \dots .$$ Here, the term with $\pm
\frac{ih}{2^n}$ appears precisely $2^n$ times. The elements
$X_i$ lie in $(\cC^\infty(S^1) \hot \cJ)^+$ since they are
exponentials of elements $ih \otimes x_i$ with $x_i \leq
2\omega \in \cJ$. This is the only place, where we use the
assumption that $\cJ$ is harmonic in an essential way.

Note that the elements $X_i$ commute with $C^{\infty}(S^1)
\otimes 1$ and with $P \otimes 1$ modulo $\cK \hot \cJ$ and
hence define abstract Kasparov $(\Cz[0,1],\cK \hot
\cJ)$-modules relative to $(D\hot \cJ)^+$. Note further that
$X_1 = W_1 \cdot W_2$ with

$$W_1 =  e^{\frac{ih}{2}} \oplus e^{-\frac{ih}{2}} \oplus 1
\oplus e^{\frac{ih}{8}} \oplus e^{\frac{ih}{8}} \oplus
e^{\frac{ih}{8}} \oplus e^{\frac{ih}{8}} \oplus
e^{-\frac{ih}{8}} \oplus \dots $$ and
$$W_2 =  e^{\frac{ih}{2}} \oplus 1 \oplus e^{-\frac{ih}{2}}
\oplus e^{\frac{ih}{8}} \oplus e^{\frac{ih}{8}} \oplus
e^{\frac{ih}{8}} \oplus e^{\frac{ih}{8}} \oplus 1 \oplus\dots
.$$

Clearly, $W_2 = (1 \otimes Y) W_1^{-1} (1 \otimes Y^{-1})$
for some permutation matrix $Y \in \cB(H)$. We see that
$[\varphi(\Cz[0,1]) \otimes 1,1 \otimes Y] = [P \otimes
1,1\otimes Y] = 0$ so that abstract Kasparov $(\Cz[0,1],\cK
\otimes \cJ)$-modules relative to $D \hot (\cB(H) \ltimes
\cJ)$ are defined, since $\cK \hot \cJ$ is a closed ideal in
$D \hot (\cB(H) \ltimes \cJ)$. (For a definition of $\cB(H)
\ltimes \cJ$ see section \ref{opideals}.) By lemma
\ref{modprop}(a) wee see that $E(\varphi \otimes 1, X_1, P
\otimes 1 ) = E(\varphi \otimes 1, W_1 1 \otimes Y W_1^{-1}(1
\otimes Y)^{-1} , P \otimes 1)=0$. A similar reasoning applies
to $X_2$. This finishes the proof. \eproof

\section{$K_0$ of a stable algebra.}
\subsection{Stabilized functors}

In this section we consider split-exact and $M_2$-stable
functors $E': \lca \to \ab$, defined on the category of
locally convex algebras.

\begin{definition}
Let $A$ be a  locally convex algebra and let $E'$ be a
$M_2$-stable and split-exact functor. Let $\cJ$ be a Banach
ideal. The algebra $A$ is called weakly $\cJ$-stable with
respect to the functor $E'$ if the functor $B \mapsto E'(B
\hot A)$ is weakly $\cJ$-stable (see definition \ref{dexact}).
\end{definition}

Proposition \ref{weaklystable} shows that $\cJ$ is weakly
$\cJ$-stable, so that there is always one obvious weakly
$\cJ$-stable algebra. Moreover, if $B$ is weakly $\cJ$-stable
and $A$ is any locally convex algebra, the $A\hot B$ is
weakly $\cJ$-stable. \mn Let now $\cJ$ be a harmonic Banach
ideal and $A$ be a weakly $\cJ$-stable algebra with respect
to the functor $E'$. In this section we show, as a corollary
of theorem \ref{maintheo}, that the associated $A$-stabilized
functor $E = E'(? \hot A): \lca \to \ab$ satisfies diffotopy
invariance. See the remark after definition \ref{m2stable}
for important examples of split-exact and $M_2$-stable
functors to which such a result could be applied.

\bprop \label{weaklystable} The functor $E: \lca \to \ab$
which assigns $A \mapsto E'(A \hot \cJ)$ is weakly
$\cJ$-stable in the sense of definition \ref{dexact}, i.e.
the algebra $\cJ$ is weakly $\cJ$-stable with respect to any
$M_2$-stable and split exact functor. \eprop

\bproof The natural map $\theta:\,\cJ(H)\hot \cJ(H) \lori
\cJ(H\hot H)$ induces a natural map

$\theta_A: E(A \hot \cJ)\to E(A)$ for every locally convex
algebra $A$. We want to show that $\theta_A \circ
j_A=\id_{E(A)} $ for the natural map $j_A: E(A) \to E(A \hot
\cJ)$ induced by the inclusion $j: \Cz \to \cJ$, i.e. $j_A =
E'(id_A \hot j \hot id_{\cJ})$.

There is an isometry $V$ in $\cB(H\hot H)$ such that $\theta
\circ (j \hot id_\cJ)(x)= \Ad (V)=VxV^*$ for $x\in \cJ $.
Choose a second isometry $V'$ in $\cB(H\hot H)$ such that
$VV^*+V'V'^*=1$. Denote by $O_2$ the algebra generated
algebraically by $1 \hot V, 1\hot V^*$ and $1 \hot V',1 \hot
V'^*$ in $A \hot \cB(H \hot H)$ and by $D$ the algebra
generated by $A \hot \cJ$ together with $1 \hot V, 1\hot V^*$
and $1 \hot V',1 \hot V'^*$ inside $A \hot \cB(H\hot H) $. We
have a split extension $$ 0 \to A \hot \cJ \to D\to O_2\to
0$$ It is easy to see that there is a unitary $U$ in
$M_2(O_2)$ and hence in $M_2(D)$ such that $UxU^* = (1 \hot
V)x(1 \hot V^*)$ for $x$ in the subalgebra
$$\left(\begin{array}{cc} A \hot \cJ & 0\\ 0 & 0
\end{array}\right)$$ of $M_2(A \hot \cJ)$. Thus, by $M_2$
stability of $E$ it is clear that $E(\Ad (1 \hot V))=\id$ on
$E(D)$. On the other hand, from the split extension above, we
see that $E(D)=E(A \hot \cJ)\oplus E(O_2)$ so that the
restriction of $\Ad (1 \hot V)$ induces the identity on $E(A
\hot \cJ)$ also. \eproof
\begin{corollary} \label{cordiff} Let $E':\lca \to \ab$ be an
$M_2$-stable and split exact functor and let $A$ be a weakly
$\cJ$-stable algebra with respect to the functor $E'$, for a
fixed harmonic Banach ideal $\cJ$. Under these circumstances,
the functor $E=E'(? \hot A)$ is diffotopy
invariant.\end{corollary} \bproof This follows now from
\ref{maintheo}.\eproof

\subsection{Algebraic $K$-theory}

We continue by applying Corollary \ref{cordiff} to the
algebraic $K$-theory functor $K_0$ and identify algebraic
$K$-theory of stable algebras with a suitable group of
homotopy classes of maps.

Consider the algebra $Q\Cz= \Cz \ast \Cz$ and denote by $e,
\bar{e}$ the two generators $e=\iota_1(1),
\bar{e}=\iota_2(1)$. The argument in \cite{CuKK}, 3.1, shows
that $E(Q\Cz\hot A)\cong E(A)\oplus E(A)$ and that $E(q\Cz\hot
A)\cong E(A)$ for each $M_2$-stable and diffotopy invariant
functor $E$. With the classical description of $K_0$ in
\cite{Bass}, the generator of $K_0(q\Cz \hot \cJ)$ is given
by the difference of equivalence classes of the idempotent
elements $p$ and $\bar{p}$ in $M_2((q\Cz \hot \cJ)^+)$:
$$
 p= W\left(\begin{array}{cc}
 \bar{e}^{\perp} & 0\\
 0 & e\end{array}\right) W^{-1} \qquad \mbox{where} \qquad
 W=\left(\begin{array}{rl}
     \bar{e}^{\perp} & \bar{e}\\
 -\bar{e} & \bar{e}^{\perp} \end{array} \right)
$$
and,
$$
\bar{p}= \left(\begin{array}{ll}
 1 & 0\\
 0 & 0 \end{array}\right).
$$
Note that $p-\bar{p} \in M_2(q\Cz \hot \cJ)$ and that
therefore $[p]-[\bar{p}] \in K_0(q\Cz \hot \cJ)$.

\bprop\label{K0} For each locally convex algebra $A$, one has
natural isomorphisms
$$K_0(A \hot \cJ) \cong \langle q\Cz, A \hot \cJ \hot M_{\infty} \rangle \cong
\lim_{\mathop{\lori}\limits_{n}} \langle q\Cz, A \hot \cJ
\hot M_n \rangle .$$ \eprop \bproof First of all, note that
the righthand side is an abelian semigroup by block sum. The
second isomorphism follows from the fact that $q\Cz$ is
finitely generated and from properties of the fine topology on
$M_\infty$.

We construct a map $\phi: \langle q\Cz, A \hot \cJ \hot
M_{\infty} \rangle \to K_0(A \hot \cJ).$ It sends a diffotopy
class $\langle h \rangle \in \langle q\Cz, A \hot \cJ \hot
M_{\infty} \rangle$ to the difference of equivalence classes
of projections $[h(p)]-[h(\bar{p})]$ in $K_0((A\hot \cJ)^+)$.
The difference lies in the direct summand $K_0(A\hot \cJ)$.
The map is well-defined by the preceding corollary.

This map is surjective since all generators of $K_0^{alg}(A
\hot \cJ)$ and their negatives are hit and since the map is
clearly a map of abelian semigroups. In order to prove
injectivity we use Lemma $7.1$ in \cite{CuDoc} which we
reprove for sake of completeness.

\begin{lemma}\label{tech}
Let $\varphi:q\Cz\to M_2(q\Cz)$ be the restriction of the
homomorphism $Q\Cz\to M_2((Q\Cz)^+)$, which sends $e$ to $p$
and $\bar{e}$ to $\bar{p}$. Then $\varphi$ is diffotopic to
the inclusion map $\iota:q\Cz\to M_2(q\Cz)$.
\end{lemma}
\bproof Let $\gamma_t:q\Cz\to M_2(q\Cz),\, t\in [0,\pi/2]$ be
the restriction of the homomorphism $\gamma'_t:Q\Cz\to
M_2((Q\Cz)^+)$ which is defined by \bglnoz
 & \gamma'_t(e)= W_t\left(\begin{array}{cc}
 \bar{e}^{\perp} & 0\\
 0 & e \end{array}\right) W_{-t}\\
 &\\
 & \gamma'_t(\bar{e}) = W_t \left(\begin{array}{cc}
 \bar{e}^{\perp} & 0\\
 0 & \bar{e} \end{array}\right) W_{-t} \eglnoz where
$$ W_t=\left(\begin{array}{cc} \bar{e}^{\perp} & 0\\ 0 &
 \bar{e}^{\perp} \end{array}\right)  +
 \left(\begin{array}{cc}
 \bar{e}\cos t & \bar{e}\sin t\\
 -\bar{e}\sin t & \bar{e}\cos t
 \end{array} \right)
 $$
For each $t$ the difference $\gamma'_t(e)-\gamma'_t(\bar{e})$
lies in the ideal $M_2(q\Cz)$. Therefore $\gamma_t$ defines a
diffotopy, which connects $\varphi$ with $\iota$. \eproof

To prove injectivity we now use Lemma \ref{tech}. Assume that
$\eta_1, \eta_2:q\Cz\to A \hot \cJ \hot M_n$ are
homomorphisms, such that $[\eta'_1(p)] = [\eta'_2(p)]$ in
$K_0(A \hot \cJ)$ (where $\eta'_i$ denotes the induced map
$M_2(q\Cz^+) \to (A \hot \cJ)^+ \hot M_{2n}$).

This means that the there exists a projector $q \in (A \hot
\cJ)^+ \hot M_k $ and an invertible element $u \in (A \hot
\cJ)^+ \hot M_{2n+k}$ such that $u (\eta'_1(p) \oplus q)
u^{-1} = \eta'_2(p) \oplus q$. This element $u$ can even be
chosen to be connected to $1$ by a differentiable family
$u_t, t\in [1,2]$, such that $1-u_t\in A \hot \cJ \hot
M_{2n+k}$ for all $t$.

Consider the homomorphisms $\zeta_1, \zeta_2:q\Cz\to A \hot
\cJ \hot M_{2n+k}$ defined as restrictions of the maps from
$Q\Cz$ that map $e$ to $\eta'_1(p) \oplus q$ and $\eta'_2(p)
\oplus q$, respectively, and $\bar{e}$ to $\bar{p} \oplus q$.
Note that $\zeta_i=M_2(\eta_i) \circ \phi$.

According to Lemma \ref{tech}, $\zeta_1= M_2(\eta_1) \circ
\varphi$ is diffotopic to $\eta_1=M_2(\eta_1) \circ \iota$
and similarly $\zeta_2$ is diffotopic to $\eta_2$. On the
other hand, the family $\zeta_t$, $t\in [1,2]$, of
homomorphisms $q\Cz\to A \hot \cJ \hot M_{2n+k}$, obtained as
restrictions of the maps from $Q\Cz$, which map $e$ to $u_t
(\eta'_1(p) \oplus q) {u_t}^{-1}$ and $\bar{e}$ to
$\eta'_1(\bar{p}) \oplus q$, defines a diffotopy  connecting
$\zeta_1$ to $\zeta_2$.
\end{proof}

\section{Determination of $kk^{\cJ}(\Cz,A)$}
\subsection{Bivariant $kk$-theory revisited}
\label{explicitkk} We now finally have to use the explicit
definition of $kk^{\alg}$. To this end we recall some
constructions
and notation from \cite{CuWeyl}.\\
Let $V$ be a complete locally convex space. Consider the
algebraic tensor algebra
$$
T_{alg}V = V\,\oplus\, V\!\!\otimes\!\! V\,\oplus\,
V^{\otimes^{3}} \oplus\,\dots
$$
with the usual product given by concatenation of tensors.
There is a canonical linear map $\sigma:V\to T_{alg}V$ mapping
$V$ into the first direct summand. We equip $T_{alg}V$ with
the locally convex topology given by the family of all
seminorms of the form $\alpha\circ \varphi$, where $\varphi$
is any homomorphism from $T_{alg}V$ into a locally convex
algebra $B$ such that $\varphi\circ\sigma$ is continuous on
$V$, and $\alpha$ is a continuous seminorm on $B$. We further
denote by $TV$ the completion of $T_{alg}V$ with
respect to this locally convex structure.\\
For any locally convex algebra $ A$ we have the natural
extension
\begin{equation}\label{UniExt} 0 \to J A \to T A \stackrel{\pi}{\to}  A \to 0.
\end{equation} Here
$\pi$ maps a tensor $x_1\otimes x_2\otimes\ldots \otimes x_n$
to $x_1x_2 \ldots x_n \in  A$ and $J A$ is defined as Ker
$\pi$. This extension is (uni)versal in the sense that, given
any extension $0 \to  I \to  E \to  B \to 0$ of a locally
convex algebra $B$, admitting a continuous linear splitting,
and any continuous homomorphism $\alpha :A\to B$, there is a
morphism of extensions \bgl\label{cla}
\begin{array}{ccccccccc} 0 &\to &  J A & \to & T A  & \to &   A
& \to & 0\\[0.1cm]
&      &  \quad\downarrow {\scriptstyle\gamma}& &
\quad\downarrow{\scriptstyle\tau} & & \quad\downarrow
{\scriptstyle\alpha} & &\\[0.1cm] 0
& \to &  I & \to &  E & \to &  B & \to & 0
\end{array}
\egl The map $\tau :T A\to  E$ is obtained by choosing a
continuous linear splitting $s: B\to  E$ in the given
extension and mapping $x_1\otimes x_2\otimes\ldots \otimes
x_n$ to $s'(x_1)s'(x_2) \ldots s'(x_n) \in E$, where
$s'\defeq s\circ \alpha$. Then $\gamma$ is the restriction of
$\tau$.\\ Choosing $0\to JB\to TB\to B\to 0$ in place of the
second extension $0\to I\to E\to B\to 0$, we see that $A \to
JA$ is a functor, i.e. any homomorphism $\alpha :A\to B$
induces a homomorphism $J(\alpha) :JA\to JB$.\\
Using the universal property of $TA$, one can associate a
\emph{classifying map} with any linearly split extension of
locally convex algebras of the form
$$0\to I\to E_1\to E_2\to \ldots \to E_n\to A\to 0$$
We consider such an extension as a complex, denoting the
arrows (boundary maps) by $\pi_i$ and we say that it is
linearly split if there is a continuous linear map $s$ of
degree -1 such that $s\pi + \pi s = \id$. Every such
splitting $s$ induces a commutative diagram of the form
$$\begin{array}{ccccccccccccccc}0&\to
&I&\to& E_1&\to &E_2&\to& \ldots& \to &E_n&\to &A&\to& 0\\&
&\uparrow& & \uparrow& &\uparrow& & \ldots&  &\uparrow&
&\uparrow& &  \\0&\to &J^nA&\to&T(J^{n-1}A)&\to &T(J^{n-2}A)
&\to& \ldots& \to &TA&\to &A&\to& 0\end{array}$$ The leftmost
vertical arrow in this diagram is the \emph{classifying map}
for this $n$-step extension. It depends on $s$ only up to
diffotopy.\\
Recall then from \cite{CuWeyl} the following definitions.
\begin{definition}\label{K-add} Let $ A$ and $ B$ be  locally convex  algebras.
For any continuous homomorphism $\varphi: A \to  B\,$, we
denote by $\langle \varphi \rangle$ the equivalence class of
$\varphi$ with respect to diffotopy and we set
$$
\langle  A , B\rangle =\{\langle \varphi \rangle |\,\varphi
\,\,\mbox{\rm is a continuous homomorphism }  A \to B\,\}
$$
Given $n\in \Zz$ we set
$$ kk^{\rm alg}_n( A ,\,  B\,) =
\lim_{\mathop{\lori}\limits_{k}} \langle J^{k-n} A ,\,
\Sigma^kB\hot\cK\, \rangle
$$
where the inductive limit is with respect to the natural maps
$$\langle J^k A ,\, \mathcal K \hat{\otimes} B(0,1)^k\,
\rangle \to \langle J^{k+1} A ,\, \mathcal K \hat{\otimes}
B(0,1)^{k+1}\, \rangle$$ mapping the diffotopy class of
$\alpha$ to the diffotopy class of $\alpha '$, where
$\alpha'$ is defined by the commutative diagram
\[
\begin{array}{ccccccccc} 0 &\to &  J^{k+1}A & \to & TJ^k A  &
\to &   J^k A & \to & 0\\[0.1cm]
&      &  \qquad\downarrow {\scriptstyle \alpha'}& &
\quad\downarrow{\scriptstyle\tau }& & \quad\downarrow
{\scriptstyle
\alpha} & &\\[0.1cm] 0 &
\to &  B'(0,1)^{k+1} & \to &  B'(0,1)^k[0,1) & \to &
B'(0,1)^k & \to & 0
\end{array}
\]
where $B' = \mc K\hot B$. \edefin Thus, by definition of
$kk^{\cJ}$, we have
$$ kk^{\rm alg}_n( A ,\,  B\,) =
\lim_{\mathop{\lori}\limits_{k}} \langle J^{k-n} A ,\,
\Sigma^kB\hot \cJ\hot\cK\, \rangle
$$
In the sequel we will also use stabilization by the finite
matrix algebras $M_n$.

\begin{proposition} \label{alg-p}
The natural maps, induced by the inclusion $M_n\to \cK$,
$$\langle J^k A, \Sigma^k B \hot \cJ \hot M_n \rangle
\to \langle J^k A, \Sigma^k B \hot \cJ \hot \cK \rangle $$
define an isomorphism $\alpha$ of abelian groups:
$$ \lim_{\mathop{\lori}\limits_{k}} \lim_{\mathop{\lori}\limits_{n}}
\langle J^{k} A ,\, \Sigma^kB\hot \cJ \hot M_n\, \rangle
\cong kk_0^{\cJ}(A,B).$$
\end{proposition}
\bproof Let $V$ be an isometry in $\cB(H)$ and $j_n :\cJ\to
M_n(\cJ)$ the inclusion map (into the upper left corner). A
standard argument shows that $j\circ \Ad V :\cJ\to \cJ\otimes
M_n$ is diffotopic to $j$.\\ Denote by $\theta: \cJ \hot \cJ
\to \cJ$ the natural tensor product of operators and by
$\varphi: \cK \to \cJ$ the natural inclusion. We claim that
the inverse $\beta$ to the map $\alpha$ above is induced by
$$j_n\circ\theta\circ
(\id_\cJ\otimes\varphi):\cJ\hot\cK\to \cJ\otimes M_n$$ The
identity $\beta\circ\alpha=1$ follows from the fact that the
composition $$\cJ\cong \cJ\otimes M_n \to \cJ \hot \cK
\stackrel{\id_{\cJ} \hot \varphi}{\lori} \cJ \hot \cJ
\stackrel{\theta}{\to} \cJ$$ is of the form $\Ad V$ as
above.\\
To show that $\alpha\circ\beta =1$, consider the inclusions
$\iota_l,\iota_r:\cK\to \cH\hot\cK$ into the left, resp.
right, factor using the standard rank 1 projector onto the
first basis vector in $H$. By \cite{CuWeyl}, 2.2.1, $\iota_l$
and $\iota_r$ are diffotopic and in fact both diffotopic to
the
natural isomorphism $\cK\to \cK\hot\cK$.\\
We want to show that the composition of the following maps
$$\cJ \hot \cK \stackrel{\id_{\cJ} \hot \varphi}{\lori} \cJ
\hot \cJ  \stackrel{\theta}{\to} \cJ \to \cJ\hot \cK$$ is
diffotopic to the identity. This composition can also be
factored as follows
$$\cJ \hot \cK \stackrel{\id_{\cJ} \hot \iota_l}{\lori} \cJ
\hot \cK \hot \cK \stackrel{\id_{\cJ} \hot \varphi \hot
\id_{\cK}}{\lori} \cJ \hot \cJ \hot \cK \stackrel{\theta \hot
\id_{\cK}}{\to} \cJ \hot \cK.$$ Replacing in this composition
$\iota_l$ by the diffotopic map $\iota_r$ we again obtain a
map of $\cJ\hot\cK\to\cJ\hot\cK$ of the form $\Ad V\otimes
\id_\cK$ which is diffotopic to id (using the fact that
$j_2:\cK\to  \cK\otimes M_2$ is diffotopic to the natural
isomorphism $\cK\cong \cK\otimes M_2$). \eproof

\subsection{Main theorem}

Our main result in this section is the following computation.
Its proof requires some preparation and is given in section
\ref{proof}.

\btheo \label{compK} For every locally convex algebra $A$ and
for every harmonic Banach ideal $\cJ$ one has
$kk_0^{\cJ}(\Cz,A)=K_0( A \hot \cJ)$. \etheo

\bremark In particular this shows that $K_0(A \hot \cL_p)$
does not depend on $p$ for $1<p<\infty$ (see also \cite{Kar},
4.1). \eremark

\begin{corollary} \label{coeff}
Let $\cJ$ be a harmonic Banach ideal. The coefficient ring
$kk^{\cJ}_*(\Cz,\Cz)$ is isomorphic to $\Zz[u,u^{-1}]$.
\end{corollary}
\bproof By properties of $kk^{alg}$ we are reduced to a
computation of $kk^{\cJ}_0(\Cz,\Cz)$ and
$kk^{\cJ}_0(\Cz,\Sigma)$. By Theorem \ref{compK} these groups
are isomorphic to the algebraic $K$-groups $K_0(\cJ)$ and
$K_0(\Sigma \cJ)$. Both algebras appearing are smooth
sub-algebras of $C^*$-algebras whose $K$-theory is
well-known, i.e.
$$K_0(\cJ) = K_0(\cL_{\infty})= \Zz$$ and
$$K_0(\Sigma \cJ) = K_0(\cC_0(\Rz,\cL_{\infty}))=0.$$ This finishes
the proof of the corollary. \eproof

\begin{corollary}
Let $W$ be the Weyl algebra, i.e the unital algebra with two
generators $x$ and $y$ satisfying the relation $xy-yx=1$
(with the fine topology). Then for every harmonic Banach ideal
$\cJ$ we have $kk^{\cJ}_0(\Cz,W)=\Zz$ and $kk^{\cJ}_1(\Cz,W)=
0$.
\end{corollary}
\bproof This follows from the result in \cite{CuWeyl}, 12.4
in combination with \ref{coeff} above. \eproof Similarly, the
$\cJ$-stable $K$-theory of many other algebras can now be
computed as an abelian group rather than just as a module
over the coefficient ring.\mn

In order to organize the notation in the following
computations, we introduce a category $H_{\cJ}$. The objects
are locally convex algebras and the morphisms between two
locally convex algebras are given by
$$[A,B] = \lim_{\mathop{\lori}\limits_{n}}\langle A, B \hot
\cJ \hot M_n \rangle.$$ Given two diffotopy classes of
continuous homomorphisms $\phi: A \rightarrow B \hot \cJ \hot
M_n$ and $\psi: B \rightarrow C \hot \cJ \hot M_m$, their
composition in $H_\cJ$ is defined to be the diffotopy class of
$$A \stackrel{\phi}{\to} B \hot \cJ \hot M_n
\stackrel{\psi \hot id}{\lori} B \hot \cJ \hot M_n \hot \cJ
\hot M_m \stackrel{id \hot \theta}{\lori} C \hot \cJ \hot
M_{nm}.$$ It is clear that composition is associative at the
level of diffotopy classes. Moreover, note that endo-functors
like $J(?)$ and $\Sigma \hot ?$ descend to endo-functors on
$H_{\cJ}$, since, for example, there is a natural map $J(A
\hot \cJ \hot M_n) \to J(A) \hot \cJ \hot M_n$.

In particular, using the new notation, we have
$$kk^{\cJ}_i(A,B) = \lim_{\mathop{\lori}\limits_{k}} [J^{k-i}
A, \Sigma^k B]$$ by proposition \ref{alg-p} and
$$K_0(A \hot \cJ) = [q\Cz, A]$$ by proposition \ref{K0}.
\bprop For all locally convex algebras $A$ and $B$ one has
$$kk^{\cJ}_i(A,B) =\lim_{\mathop{\lori}\limits_{k}} [J^{2k-i} A,B] $$
\eprop This identity was noted in \cite{CuWeyl}, remark 8.4.
In order to give an explicit proof, we have to introduce some
notation. Throughout, we are working in the category
$H_{\cJ}$. Denote by $\rho_A : JA \to \Sigma A$ the
classifying map of the cone extension of $A$. Denote by
$\ve'_A: J\Sigma A \to A$ the classifying map of the Toeplitz
extension tensored with $A$. We define $\ve_A = \ve'_A \circ
J(\rho_A) : J^2A \to A$. (The map $\ve_A$ has the important
interpretation as the classifying map for the $2$-step
extension given by the Yoneda product of the Toeplitz
extension and the cone extension.)

Furthermore, note that $$J^2(\ve_{A}) = \ve_{J^2(A)}$$ by
Korollar 3.1.1 in \cite{CuDoc}, or by Lemma 4.6 in
\cite{CuWeyl}, for all algebras $A$. For $\phi: A \to B$, we
also have the identity
$$\phi \circ \ve_{A} = \ve_{B} \circ J^2(\phi).$$

We define inductively $\ve'^n_A: J^n \Sigma^nA \to A$ by
setting $\ve'^n_A = \ve'^{n-1}_A \circ
J^{n-1}(\ve'_{\Sigma^{n-1}A})$ and $\rho^n_A: J^n A \to
\Sigma^n A$ by setting $\rho^n_A = \rho_{\Sigma^{n-1} A}
\circ J(\rho^{n-1}_A)$. We also define $\ve^n_A: J^{2n} A \to
A$ by setting $\ve^n_A = \ve^{n-1}_A \circ \ve_{J^{2n-2}A}$.
Note that, in the definition of $\ve^n_A$, all other choices
lead to the same definition of $\ve^n_A$, since $J^2(\ve_A) =
\ve_{J^2A}$.

Note that $\ve'^n_A \circ J^n(\rho^n_A) = \ve^n_A :J^{2n}A
\to A$, as the following induction argument shows.
\[ \begin{array}{rcl}
\ve'^n_A \circ J^n(\rho^n_A) & = & \ve'^{n-1}_A \circ
J^{n-1}(\ve'_{\Sigma^{n-1}A}) \circ
J^n( \rho_{\Sigma^{n-1} A} \circ J(\rho^{n-1}_A)) \\
&=& \ve'^{n-1}_A \circ J^{n-1}(\ve'_{\Sigma^{n-1}A} \circ J(\rho_{\Sigma^{n-1} A}) \circ J^2(\rho^{n-1}_A)) \\
&=& \ve'^{n-1}_A \circ J^{n-1}(\ve_{\Sigma^{n-1}A} \circ J^2(\rho^{n-1}_A)) \\
&=& \ve'^{n-1}_A \circ J^{n-1}(\rho^{n-1}_A \circ \ve_{J^{n-1}A}) \\
&=& \ve'^{n-1}_A \circ J^{n-1}(\rho^{n-1}_A) \circ \ve_{J^{2n-2}A} \\
\end{array}\]

With this notation, the abelian group $kk^\cJ(A,B)$ is defined
as the direct limit of a system of abelian groups $[J^n A,
\Sigma^nB]$ via a stabilization map $$[J^n A, \Sigma^n B] \ni
\psi \mapsto \rho_{\Sigma^n B} \circ J(\psi) \in [J^{n+1},
\Sigma^{n+1} B].$$ The $k$-th stabilization map is given by
$$[J^nA, \Sigma^n B] \ni \phi \mapsto \rho^k_{\Sigma^n B}
\circ J^k(\phi) \in [J^{n+k} A, \Sigma^{n+k} B],$$ which,
again, follows from an easy induction argument.

The right hand side is defined as the direct limit of abelian
semi-groups $[J^{2n}A,B]$ via the stabilization map
$$[J^{2n}A,B] \in \psi \mapsto \ve_B \circ J^2(\psi)  = \psi \circ \ve_{J^{2n} A} \in [J^{2n+2}A,B].$$

\bproof First of all, we define a map $[J^n A,\Sigma^n B] \to
[J^{2n} A, B]$ by $\phi \mapsto \ve'^n_B \circ J^n(\phi) \in
[J^{2n},B].$ This assignment induces a map of directed
systems since it is compatible with the stabilization maps.
Indeed, $\rho_{\Sigma^n B} \circ J(\phi)$ is mapped to
$\ve'^{n+1}_B \circ J^{n+1}( \rho_{\Sigma^n B} \circ
J(\phi))$ and the following computation shows that this is
the desired result.
\[
\begin{array}{rcl}
\ve'^{n+1}_B \circ J^{n+1}( \rho_{\Sigma^n B} \circ J(\phi))
&=& \ve'^n_B \circ J^n( \ve'_{\Sigma^n B})  \circ J^n(J(\rho_{\Sigma^n B} \circ J^2(\phi))) \\
&=& \ve'^n_B \circ J^n( \ve'_{\Sigma^n B} \circ J(\rho_{\Sigma^n B}) \circ J^2(\phi)) \\
&=& \ve'^n_B \circ J^n(\ve_{\Sigma^n B} \circ J^2(\phi)) \\
&=& \ve'^n_B \circ J^n(\phi \circ \ve_{J^n A})\\
&=& \ve'^n_B \circ J^n(\phi) \circ \ve_{J^{2n} A}\\
\end{array} \]

This implies that there is a well-defined map
$$\sigma_{A,B}: \lim_{\mathop{\lori}\limits_{k}} [J^k A, \Sigma^k B] \to \lim_{\mathop{\lori}\limits_{k}} [J^{2k}A,B].$$

Consider a class $[\phi] \in
\mathop{\lim}\limits_{\mathop{\lori}\limits_{k}} [J^k A,
\Sigma^k B]$ which is represented by a homomorphism $\phi:
J^n A \to \Sigma^n B$. If the composition $\ve'^n_B \circ
J^n(\phi) \circ \ve^m_{J^{2n} A}: J^{2m+2n} A\to B$ is
homotopic to zero for some $m \in \Nz$, then also $[\ve'^n_B]
\circ [J^n(\phi)] \circ [\ve^m_{J^{2n} A}]=0$ as element in $
\mathop{\lim}\limits_{\mathop{\lori}\limits_{k}}[J^k
J^{2m+2n} A, \Sigma^k B].$

Note that, for any $A$ and $B$, the functor $J$ maps
$\mathop{\lim}\limits_{\mathop{\lori}\limits_{k}} [J^k A,
\Sigma^k B]$ isomorphically onto
$\mathop{\lim}\limits_{\mathop{\lori}\limits_{k}} [J^k JA,
\Sigma^k JB]$. Furthermore, $[\ve'_A]$ and $[\ve_A]$ are
invertible for all $A$. This together implies that $[\phi]=0$
as a class in
$\mathop{\lim}\limits_{\mathop{\lori}\limits_{k}} [J^{k+n} A,
\Sigma^{k+n} B]$ and hence that $\sigma_{A,B}$ is injective.

We still have to show that $\sigma_{A,B}$ is surjective.
Consider a class $[\phi] \in
\mathop{\lim}\limits_{\mathop{\lori}\limits_{k}} [J^{2k}A,B]$
which is represented by a homomorphism $\phi: J^{2n} A \to
B$. The element $\phi$ gives rise to a class in
$\mathop{\lim}\limits_{\mathop{\lori}\limits_{k}} [J^k J^{2n}
A, \Sigma^k B]$. The natural map
$$\lim_{\mathop{\lori}\limits_{k}} [J^k A, \Sigma^k B] \to
\lim_{\mathop{\lori}\limits_{k}} [J^k J^{2n} A, \Sigma^k B]$$
which is induced by the assignment
$$[J^k A ,\Sigma^k B] \ni \psi \mapsto \psi \circ J^k (\ve^n_A) \in [J^k J^{2n} A, \Sigma^k B]$$ is an isomorphism.
I.e. there is some $k \in \Nz$ and $\eta: J^k A \to \Sigma^k
B$ such that
$$ \eta \circ J^k(\ve^n_A) = \rho^k_B \circ J^k(\phi).$$
We claim that $[\eta] \in
\mathop{\lim}\limits_{\mathop{\lori}\limits_{k}} [J^k A,
\Sigma^k B]$ does the job. Indeed, it is mapped to

$\ve'^k_B \circ J^k(\eta): J^{2k} A \to B.$ Stabilizing yields
\[
\begin{array}{rcl}
\ve'^k_B \circ J^k(\eta) \circ \ve^n_{J^{2k} A}
&=& \ve'^k_B \circ J^k(\eta) \circ J^{2k}(\ve^n_A) \\
&=& \ve'^k_B \circ J^k(\rho^k_{B}) \circ J^{2k}(\phi)\\
&=& \ve^k_B \circ J^{2k}(\phi) \\
&=& \phi \circ \ve^k_{J^{2k}A}.\\
\end{array} \]
This shows the surjectivity of $\sigma_{A,B}$. It is another
easy check to show, that the assignment is compatible with
the several composition products. \eproof

\subsection{Proof of the main theorem}
\label{proof} We now proceed by proving Theorem \ref{compK}.
By proposition \ref{K0} it suffices to shows that there is a
natural isomorphism
$$[q\Cz,A] \cong \lim_{\mathop{\lori}\limits_{k}} [J^{2k}\Cz,A].$$ The existence of a natural map will become apparent in the sequel of the proof.

\bproof  There is a natural map $\alpha_1: J^2(\Cz) \to q\Cz$
given as the classifying map of a $2$-step extension which
comes as the Yoneda-product of the Toeplitz extension and a
canonical extension

$$0 \to \Sigma q\Cz \to E \to \Cz \to 0,$$
where $E = \{f: [0,1] \to Q\Cz| f(0) \in \Cz \ast 0, f(1) \in
0\ast \Cz, f(0) - f(t) \in q\Cz, \forall t\in[0,1]\}.$

Denote by $\delta: q\Cz \to \Cz$ the restriction of $id \ast
0: \Cz \ast \Cz \to \Cz$. It follows from the definition of
$\ve_{\Cz}: J^2(\Cz) \to \Cz$ that $\delta \circ \alpha_1 =
\ve_{\Cz}$.

We define natural maps $\alpha_n = \alpha_1 \circ
\ve^{n-1}_{J^2\Cz} \in [J^{2n}\Cz,q\Cz]$. \mn Since $A
\mapsto K_0(A \otimes \cJ)$ is a half-exact and diffotopy
invariant functor by Theorem \ref{maintheo}, we conclude by
Lemma \ref{halfexact} that there are natural isomorphisms
$K_0(J^{2n}\Cz \hot \cJ) \cong K_0(\Sigma^{2n} \cJ)$.
However, the righthand side of the last equation is
isomorphic to $K_0(\cC_0(\Rz^{2n}))= \Zz$ and contains a
canonical generator. Denote the canonical generators of
$K_0(J^{2n}\Cz \hot \cJ)$ by $\beta_m \in [q\Cz,
J^{2n}(\Cz)]$. The following identities are immediate, once
we have the alternative description of $K_0(?\hot \cJ)$,
given in proposition \ref{K0}.

\[
\begin{array}{rcl}
\alpha_n \circ \beta_n & = & id_{q\Cz} \\
\ve^k_{J^{n-k}\Cz} \circ \beta_n & = & \beta_{n-k} \\
\ve^n_{\Cz} \circ \beta_n & = & \delta
\end{array}\]

It is clear from the identities above that the maps $\alpha^*
: {\rm lim}_{n\in \Nz} [J^{2n}\Cz,A] \to [q\Cz,A]$ and
$\beta^* : [q\Cz,A] \to {\rm lim}_{n\in \Nz} [J^{2n}\Cz,A]$,
which are induced from $\alpha_n$ and $\beta_n$ by
precomposition, are well defined. Furthermore, $\beta^* \circ
\alpha^*$ is equal to the identity. We now show that
$\alpha^* \circ \beta^*$ is also equal to the identity.
\[
\begin{array}{rcl}
\beta_n \circ \alpha_n \circ \ve^n_{J^{2n}\Cz} & = & \ve^n_{J^{2n}\Cz} \circ J^{2n}(\beta_n \circ \alpha_n) \\
& = & J^{2n}(\ve^n_\Cz) \circ J^{2n}(\beta_n \circ \alpha_n) \\
& = & J^{2n}(\ve^n_\Cz \circ \beta_n \circ \alpha_n) \\
& = & J^{2n}(\delta \circ \alpha_n) \\
& = & J^{2n}(\delta \circ \alpha_1 \circ \ve^{n-1}_{J^2\Cz}) \\
& = & J^{2n}(\ve^n_{\Cz}) \\
& = & \ve^n_{J^{2n}\Cz}
\end{array}\]
\eproof

\section{Computation of $K_1(\cJ)$ for harmonic Banach ideals}

In this section we want to give a computation of the
algebraic $K$-theory group $K_1(\cJ)$ for a harmonic Banach
ideal $\cJ$ (see definition \ref{har}). The result
contradicts an old result in \cite{Kar} prop. 4.1. The error
in the proof of proposition 4.1 in \cite{Kar} was brought to
our attention by Valqui and Corti\~{n}as. This concrete
computation fits nicely with the far more general structure
theorem about higher algebraic $K$-theory of locally convex
algebras stabilized by harmonic Banach ideals, which was
obtained by Corti\~{n}as in \cite{Cor} using our diffotopy
invariance theorem for weakly $\cJ$-stable, $M_2$-stable,
split-exact functors. \mn Let $A$ be a locally convex
algebra. In this section $A^+$ denotes the unitization by
$\Zz$ rather than by $\Cz$. A priori, this difference matters
and one has to be careful not to use the complex unitization.
\mn The abelian group $K_1(\cJ)$ is defined as $\ker
(K_1(\cJ^+) \to K_1(\Zz))$. Denote by $\cJ^2$ the (algebraic)
square of $\cJ$, i.e. the image of the algebraic tensor
product $\cJ\otimes\cJ$ under the multiplication map. We
consider the map $\cJ^+ \to \cJ^+/{\cJ^2}$. Since
$\cJ^+/{\cJ^2}$ is abelian, we have a naturally defined
determinant map
$$\bigcup_{n\in \Nz} Gl_n(\cJ^+/\cJ^2) \to (\cJ^+/\cJ^2)^{\times} = \Zz^\times \ltimes \cJ/\cJ^2.$$ We get induced maps
$$det: K_1(\cJ) \to K_1(\cJ/\cJ^2) \to \ker(\Zz^\times \ltimes \cJ/\cJ^2 \to \Zz^\times)
= \cJ/\cJ^2$$ (we use here the identification of the additive
group $\cJ/\cJ^2$ with the multiplicative group $\{1+a\;| \;
a\in \cJ/\cJ^2\}$).\mn We want to show that $det: K_1(\cJ) \to
\cJ/\cJ^2$ is an isomorphism. Note, that it is obviously
surjective, since $det\,(e^a)= a +\cJ^2$ for any $a \in \cJ$.

\begin{stheorem} Let $\cJ$ be a harmonic Banach ideal.
The natural determinant map yields an isomorphism
$$K_1(\cJ) \cong \cJ/{\cJ^2}.$$
\end{stheorem}

The preceding result implies that our diffotopy result does
not extend to higher algebraic $K$-theory. The first
topological $K$-theory of the Schatten ideals is well-known
to be zero. The proof of the preceding theorem is given after
stating and proving a lemma.

Note that, every class $[z] \in K_1(\cJ)$ is represented by
an element $z \in Gl_n(\cJ^+)$ which maps to the identity
under the canonical evaluation onto $Gl_n(\Zz)$. Using the
next lemma, we are able to show injectivity of the
determinant map.

\begin{slemma}
\begin{itemize}
\item[(a)]
Every invertible element in $M_n(\cJ)^+$ which maps to $1 \in
\Zz$ is a product of exponentials of elements in $M_n(\cJ)$.

An invertible element in $M_n(\cJ^2)^+$ which is connected by
a norm continuous path of elements, which are invertible in
$M_n(\cJ^2)^+$, to $1 \in M_n(\cJ^2)^+$ is a product of
exponentials of elements in $M_n(\cJ^2)$.
\item[(b)]
Let $g \in M_n(\cJ^2)$. The element $e^g$ is invertible in
$M_n(\cJ)^+$ and $[e^g] = 0$ in $K_1(\cJ)$.
\item[(c)]
Let $g,h \in M_n(\cJ)$. We have that $[e^{g+h}] = [e^g] +
[e^h]$ in $K_1(\cJ)$.
\item[(d)]
Let $g \in M_n(\cJ)$. If $tr(g) \in \cJ^2$, then $[e^g]=0$ in
$K_1(\cJ)$.
\end{itemize}

\end{slemma}
\bproof Since $M_n(\cJ)$ is isomorphic to $\cJ$ and
$(M_n\cJ)^2 = M_n \cJ^2$  we can restrict our reasoning, for
a proof of (a),(b) and (c), to the case $n=1$. \mn (a) If $a
\in \cJ^+$ is close to one, then the logarithmic series
converges to an element in $\cJ$. An easy calculation shows
that the logarithm lies in $\cJ^2$ if $a \in (\cJ^2)^+$ (it
can be written as a product of $a$ by an element in $\cJ^+$).
The assertion follows by standard arguments using compactness
and the fact that the group of invertible elements in $\cJ^+$
which map to $1 \in \Zz$ is connected. \mn (b) The proof of
this lemma follows the idea in \cite{Kar}, 4.1, using a
scheme going back to \cite{PT}. As in the proof of the
homotopy invariance theorem we consider certain invertible
elements which we want to represent by commutators. Consider
$M_2(\cJ)$ as acting on $H \oplus H$ and choose an isometry
between the second copy of $H$ in this direct sum and
$\oplus_{n\in \Nz}H$. Under this identification, we define
$$X_1=e^g \oplus e^{-g/2} \oplus e^{-g/2} \oplus e^{g/4}
\oplus \dots$$ and
$$X_2=1 \oplus e^{g/2} \oplus e^{g/2} \oplus e^{-g/4} \oplus \dots $$
(again, with $e^{\pm g/2^n}$-term appearing $2^n$ times) can
be considered as elements in $M_{2}(\cJ)$, since $\cJ$ is
harmonic (compare to the proof of Theorem \ref{maintheo}).
Clearly $X_1 X_2 = e^g \oplus 1$. In order to show that
$[e^g]=0$, it suffices to show $[e^g\oplus 1]=0$ in
$K_1(\cJ)$. We want to show that $X_1$ and $X_2$ are products
of commutators. This implies that $e^g \oplus 1$ is also a
product of commutators and hence finishes the proof.

We first concentrate on $X_1$. As before we construct matrices
$$W_1 = e^{g/2} \oplus e^{-g/2} \oplus 1 \oplus e^{g/8} \dots $$ and
$$W_2 = e^{g/2} \oplus 1 \oplus e^{-g/2} \oplus e^{g/8} \dots $$
which satisfy $W_1 W_2 = X_1$. There is an explicit
isomorphism $M_{2}(\cJ)^+ \cong M_3(\cJ)^+$ (preserving the
ideal $\cJ^2$) such that the matrix $W_1$ is mapped to a
matrix of the form $\gamma \oplus \gamma^{-1} \oplus 1$. The
element $h=1-\gamma$ is in $\cJ^2 $ (since the exponential
series was already convergent in $(\cJ^2)^+$) and therefore
decomposes into $\alpha \beta$ with $\alpha, \beta \in \cJ$.
To see this assume that $h = \sum_{i\leq n} a_i b_i$ with
$a_i$ and $b_i$ in $\cJ$. Then there are $A$ and $B$ in
$M_n\cJ$ such that $h \oplus 0 = AB$. Identify then $\cJ$ with
$M_n(\cJ)$ using $\Ad V$ for a suitable isometry $V$.\\
By the proof of Vaserstein's lemma (or by direct
computation), we have that

$$ \gamma \oplus \gamma^{-1} \oplus 1 = 
{\scriptstyle \left [ \scriptstyle \left( 
\begin{array}{ccc} \scriptstyle \gamma & \scriptstyle0 & \scriptstyle0 \\ \scriptstyle 0 & \scriptstyle 1 & \scriptstyle 0
\\ \scriptstyle -\gamma \beta & \scriptstyle 0 & \scriptstyle 1 \end{array} \scriptstyle \right) , \scriptstyle \left(
\begin{array}{ccc} \scriptstyle \gamma^{-1} & \scriptstyle 0 & \scriptstyle \gamma^{-1} \scriptstyle \alpha \\ \scriptstyle 0
& \scriptstyle 1 & \scriptstyle 0 \\ \scriptstyle 0 & \scriptstyle 0 & \scriptstyle 1 \end{array} \scriptstyle \right) \scriptstyle \right] 
\scriptstyle \left[ \scriptstyle \left( \begin{array}{ccc} \scriptstyle 1 & \scriptstyle 0 & \scriptstyle 0  \\ \scriptstyle 0 & 
\scriptstyle \gamma^{-1} &\scriptstyle -\gamma^{-1} \scriptstyle \alpha \\ \scriptstyle 0 & \scriptstyle 0 & \scriptstyle 1 \end{array} \scriptstyle \right),
\scriptstyle \left( \begin{array}{ccc} \scriptstyle 1 & \scriptstyle 0 & \scriptstyle 0 \\ \scriptstyle 0 & \scriptstyle \gamma & 
\scriptstyle 0 \\ \scriptstyle 0 & \scriptstyle \gamma \scriptstyle \beta & \scriptstyle 1 \end{array} \scriptstyle \right) \scriptstyle \right]} $$

where $[\cdot,\cdot]$ denotes a multiplicative commutator.
Thus, $W_1$ is a product of two multiplicative commutators in
$M_2(\cJ)^+$. A similar reasoning applies to $W_2$ and
matrices occurring in a similar decomposition of $X_2$. Thus
$[e^g]=0$ in $K_1(\cJ)$.\mn (c) Since the leading terms in the
series expansion of $e^{g+h}e^{-g}e^{-h}-1$ vanish, this
element can be written as a sum of four terms where each term
is a product of $g^2, gh, hg, h^2$ respectively, by an
element in $\cJ^+$. Therefore $e^{g+h}e^{-g}e^{-h}$ as well
as its inverse lie in $(\cJ^2)^+$. By (a)
$e^{g+h}e^{-g}e^{-h}$ is a product of exponentials of
elements in $(\cJ^2)^+$. Using (b), this implies the claim.
\mn (d) A matrix with trace in $\cJ^2$ is a finite sum of
\begin{itemize}\item off-diagonal matrices with one entry in
$\cJ$,\item matrices of the form $a \otimes e_{ii} - a
\otimes e_{jj}$ with $a \in \cJ$ and
\item a matrix $a \otimes e_{11}$ with $a\in
\cJ^2$.\end{itemize} The classes in $K_1(\cJ)$ of their
exponentials are zero, whence the claim by iterated
application of (c). \eproof

We now proceed with the proof of the theorem.

\bproof Let $A$ be an invertible element in $M_n(\cJ^+)$
which maps to the identity in $M_n(\Zz)$ and with determinant
zero. By (a) of the preceding lemma it is of the form $e^{h_1}
\cdots e^{h_k}$ and $det\, A= h_1+ \dots +h_k + \cJ^2 =
\cJ^2$. The matrix $e^{h_1} \oplus \cdots \oplus e^{h_n} \in
M_{nk}(\cJ)^+$ has the same class in $K_1(\cJ)$ and the trace
of its logarithm is just a lift of its determinant $det\, A$
to $\cJ$ and hence in $\cJ^2$. By (d) of the preceding lemma
the class in $K_1(\cJ)$ is zero. This shows injectivity of the
determinant map. Surjectivity was obvious, hence the
assertion. \eproof

\bremark It is of course always true that $K_1(\cJ)$ maps
surjectively onto $\cJ/[\cJ,\cJ]$, $\cJ$ being a Banach
ideal. \eremark

\end{document}